\documentclass[a4paper,11pt,
 leqno
]{article}

\usepackage[utf8]{inputenc}
\usepackage{amsthm}
\usepackage{amssymb}
\usepackage{mathrsfs}
\usepackage{amsmath}
\usepackage{mathtools}
\usepackage{enumitem}
\usepackage{bbm}
\setlist{topsep=3pt,itemsep=1pt,parsep=0pt}
\usepackage[a4paper]{geometry}
\usepackage[obeyDraft]{todonotes}
\usepackage{centernot}
\usepackage{hyperref}
\hypersetup{linktocpage,
  colorlinks   = true,
  urlcolor     = blue,
  linkcolor    = blue,
  citecolor   = blue
}\usepackage[all]{hypcap}

\usepackage{titlesec}
\titleformat{\subsection}[runin]{\normalfont\bfseries}{\thesubsection.}{.5em}{}[.]
\titlespacing{\subsection}{0pt}{2ex plus .1ex minus .2ex}{.8em}

\title{
\normalsize
\textbf{ON THE ONE-ARM PROBABILITY FOR THE METRIC GRAPH\\ GAUSSIAN FREE FIELD IN HIGH DIMENSIONS}}
\author{}
\date{}
\numberwithin{equation}{section}
\allowdisplaybreaks[1]

\let\ophi\phi
\newcommand{\Ran}{\mathcal{I}}
\newcommand{\rL}{L}
\newcommand{\triv}{\textnormal{triv}}
\newcommand{\nt}{\textnormal{nt}}
\newcommand{\ran}{\textnormal{ran}}

\newcommand{\Z}{\mathbb Z}
\newcommand{\R}{\mathbb R}
\renewcommand{\P}{\mathbb{P}}
\newcommand{\E}{\mathbb{E}}
\renewcommand{\epsilon}{\varepsilon}
\renewcommand{\tilde}{\widetilde}

\renewcommand{\phi}{\varphi}
\newcommand{\capacity}{\textnormal{cap}}
\newcommand{\Law}{\operatorname{Law}}
\newcommand{\cC}{\mathscr{C}}
\newcommand{\cG}{\mathcal{G}}
\newcommand{\defeq}{\stackrel{\text{def.}}{=}}
\makeatletter
\newcommand{\nlr}{\mathrel{\mathpalette\nlr@\relax}}
\newcommand{\nlr@}[2]{\ooalign{$\m@th#1\leftrightarrow$\cr\hidewidth\raisebox{.2ex}{\scalebox{1}[.75]{$\m@th\nlr@sub#1/$}}\hidewidth\cr}}
\newcommand{\nlr@sub}[1]{\ifx#1\displaystyle\scriptstyle\else\ifx#1\textstyle\scriptstyle\else\scriptscriptstyle\fi\fi}
\makeatother

\newtheorem{theorem}{Theorem}[section]
\newtheorem{lemma}[theorem]{Lemma}
\newtheorem{proposition}[theorem]{Proposition}
\newtheorem{corollary}[theorem]{Corollary}
\theoremstyle{definition}

\newtheorem{remark}[theorem]{Remark}

\newcounter{constant}
\makeatletter
\newcommand{\CL}[1]{\stepcounter{constant}\expandafter\xdef\csname cst@#1\endcsname{\theconstant}C_{\theconstant}}
\newcommand{\CR}[1]{C_{\csname cst@#1\endcsname}}
\makeatother

\usepackage{color}
\definecolor{Red}{rgb}{1,0,0}
\definecolor{Blue}{rgb}{0,0,1}
\definecolor{Olive}{rgb}{0.41,0.55,0.13}
\definecolor{Yarok}{rgb}{0,0.5,0}
\definecolor{Green}{rgb}{0,1,0}
\definecolor{MGreen}{rgb}{0,0.8,0}
\definecolor{DGreen}{rgb}{0,0.55,0}
\definecolor{Yellow}{rgb}{1,1,0}
\definecolor{Cyan}{rgb}{0,1,1}
\definecolor{Magenta}{rgb}{1,0,1}
\definecolor{Orange}{rgb}{1,.5,0}
\definecolor{Violet}{rgb}{.5,0,.5}
\definecolor{Purple}{rgb}{.75,0,.25}
\definecolor{Brown}{rgb}{.75,.5,.25}
\definecolor{Grey}{rgb}{.7,.7,.7}
\definecolor{Black}{rgb}{0,0,0}

\def\black{\color{Black}}

\begin{document}
\maketitle
\vspace{0.1cm}
\begin{center}
\vspace{-1.9cm}
Pierre-Fran\c{c}ois Rodriguez\footnotemark[1]$^{,}$\footnotemark[2], Jialu Shi\footnotemark[1] and Wen Zhang\footnotemark[2]
\footnotetext[1]{University of Cambridge, Department of Pure Mathematics and Mathematical Statistics, Wilberforce Road, Cambridge CB3 0WB, UK, \url{pfr26@cam.ac.uk}, \url{js3007@cam.ac.uk}}
\footnotetext[2]{Imperial College London, Department of Mathematics, London SW7 2AZ, UK, \url{kz2623@ic.ac.uk}}

\end{center}
\vspace{0.1cm}
\begin{abstract}
\centering
\begin{minipage}{0.80\textwidth}
We consider the sign clusters of the Gaussian free field on the cable system of a transient weighted graph with volume growth exponent $\alpha$ and Green function decay exponent $\nu$. We prove that the probability that the cluster of a vertex reaches distance $R$ is of order $R^{-\beta}$ when $2\alpha<3\nu$, where $\beta=\alpha-\nu$ is the walk dimension, and that the same holds up to a logarithmic factor in the upper bound when $2\alpha=3\nu$. On $\mathbb Z^d$ these conditions read $d>6$ and $d=6$. The upper bound combines a general entropic bound, which compares the law of the cluster under the critical loop soup and under a massive one, with a Simon--Lieb-type inequality for the massive field derived from its Markov property. The arguments are robust and require no transitivity. The lower bound follows from a competition between the capacity and the volume of the cluster. In both bounds, the key quantity is the expected exit time of the random walk from balls, which links the volume of the cluster to its radius.
\end{minipage}
\end{abstract}
\thispagestyle{empty}

\vfill
\begin{flushright}
\today
\bigskip
\end{flushright}

\newpage
\setcounter{footnote}{0}

\section{Introduction}\label{sec:intro}

 Converting point-to-point estimates into point-to-set ones is a classical problem in percolation. In \cite{kozma_nachmias_2011}, Kozma and Nachmias showed that for Bernoulli percolation on $\Z^d$ in dimension $d>6$ the radius tail (the so-called `one-arm' probability) decays like $R^{-2}$, given two-sided bounds on the critical two-point function; see \cite{asselah_schapira_sousi_2025,vEGPS,panis_schapira_2026} for recent related developments. For the sign clusters of the Gaussian free field (GFF) on the cable system of $\Z^d$, the analogous result was obtained by Cai and Ding \cite{cai_one-arm_2025} by adapting the strategy of \cite{kozma_nachmias_2011}. At the heart of most of these arguments lies a delicate regularity analysis of the cluster, whose purpose is to relate its intrinsic cluster geometry to that of the ambient lattice.

The GFF comes with its own set of obstructions (which the authors of \cite{cai_one-arm_2025} skillfully overcame), but also with its own bag of tools. In this article we give a short proof of the one-arm exponent for this model on general transient graphs, with no symmetry assumption, and up to a logarithm in the critical dimension. The regularity analysis disappears altogether: this is essentially due to the Markov property of the field, which allows (with some work, involving a very natural perturbation of the cluster using a small mass $\kappa > 0$) to express the probability that the cluster exits a ball $D$ directly in terms of the random walk $X$ on $\cG$. The `pioneers' (i.e. the boundary points of a region connected to the origin inside it) implicit in \cite{kozma_nachmias_2011}, see also \cite{duminil-copin_tassion_2016, panis_schapira_2026}, are thus, so to speak, generated by the exit time $T_D$ (and location) of $X$ from $D$, which feature prominently in both upper and lower bounds derived below (see \eqref{eq:lb} and \eqref{eq:ub}), and drive the contraction mechanism producing sharp bounds on the one-arm probability.

We now describe our results. We consider a transient weighted graph $\cG=(G,\lambda)$, where $G$ is a countably infinite set and
$\lambda_{x,y}=\lambda_{y,x}\geq0$, $\lambda_{x,x}=0$, $x,y\in G$, are weights such that the induced
graph with edge set $E=\{\{x,y\}:\lambda_{x,y}>0\}$ is connected and locally finite. Throughout, we always make the natural, see \cite[Definition~1.2]{Barlow_2017} (and mild) ellipticity assumption
\begin{align}
&\tag{$p_0$} p_{x,y}\geq c_0 \text{ for all } x\sim y \label{eq:p_0}.
\end{align}
We write $x\sim
y$ if $\{x,y\}\in E$, $\lambda_x\coloneqq\sum_{y}\lambda_{x,y}$, $d_{G}$ for the graph distance, and
for $D\subset G$, $D^c\coloneqq G\setminus D$ and $\partial D\coloneqq\{y\in D^c: y\sim z \text{ for
some } z\in D\}$. Let $X$ be the continuous-time Markov chain on $G$ jumping from $x$ to $y$ at rate 
$p_{x,y}\coloneqq\lambda_{x,y}/\lambda_x$, with law $P_x$ when $X_0=x$, and
\begin{equation}\label{eq:green}
p_t(x,y)\defeq\frac{P_x(X_t=y)}{\lambda_y}, \qquad g_{\kappa}(x,y)\defeq\int_0^\infty e^{-\kappa t} p_t(x,y)\,dt, \qquad
x,y\in G,\ t\geq0,\ \kappa \geq0.
\end{equation}
We abbreviate $g=g_0$. For finite $K\subset G$, we let $\capacity(K)\coloneqq\sum_{x\in K}\lambda_xP_x(\tilde H_K=\infty)$, where
$\tilde H_K$ is the return time to $K$ of the discrete skeleton of $X$. Let $\tilde\cG$ be the metric graph obtained by replacing each edge $\{x,y\}\in E$ by an interval of
length $\rho_{x,y}\coloneqq(2\lambda_{x,y})^{-1}$, and $\tilde\varphi$ the Gaussian free field on
$\tilde\cG$ with law $\P$, whose covariance restricted to $G$ is $g$. Let
$\cC_x$ denote the cluster of $x\in G$ in $\{\tilde\varphi\geq0\}$ ($\cC_x=\emptyset$ if
$\tilde\varphi_x<0$), write $x\leftrightarrow A$ if $\cC_x\cap A\neq\emptyset$ for $A\subset G$.

 For a (finite) set $D\subset G$, let $T_D\coloneqq\inf\{t\geq0:X_t\notin D\}$ denote the exit time of $D$ by $X$. The expected exit time $E_x[T_D]$ will play an important role in this article (it will soon re-appear in Proposition~\ref{prop:ub}, below). In the sequel, $c,c',C,\dots$ denote positive constants which may
change from place to place and depend only on the constants in the conditions in force.

\begin{proposition}[Lower bound]\label{prop:lb}
Assume \eqref{eq:p_0} and $(\lambda_x\vee g(x,x))\leq C$ for all $x\in G$. Then for all $x\in G$ and finite $D\ni x$,
\begin{equation}\label{eq:lb}
\P\big(x\leftrightarrow D^c\big)\geq c\max\big\{\capacity(D)^{-1/2}, E_x[T_D]^{-1}\big\}.
\end{equation}
\end{proposition}

 The assumptions of Proposition~\ref{prop:lb} on $\mathcal G$ will appear frequently below. Note that, up to modifying the constants, the condition $(\lambda_x\vee g(x,x))\leq C$ is equivalent to the \emph{two-sided} bounds $c\leq\lambda_x\leq C$ and $c\leq g(x,x)\leq C$: one direction is obvious, and for the converse one uses that $g(x,x)\geq\lambda_x^{-1}$, which follows by considering the time spent at $x$ before the first jump of $X$.

The proof of Proposition~\ref{prop:lb} is straightforward; see Section~\ref{sec:prelim}. It is based on available information on the cluster capacity and two--point functions of the model. A variant of \eqref{eq:lb} under weaker assumptions on $\cG$ can be found in Remark~\ref{rem:lb_general}. As discussed there, \eqref{eq:lb} entails a competition between capacity and volume of the cluster in $D$.

To get a feeling for the meaning of \eqref{eq:lb}, and for later reference, we consider the following conditions on $\mathcal G$, which allow to draw conclusions about $\capacity(D)$ and $E_x[T_D]$ when $D$ is a ball: for a distance function
$d$ on $G$ with balls $B(x,R)\coloneqq\{y\in G:d(x,y)\leq R\}$, and two parameters $\alpha>2$,
$\nu\in(0,\alpha-2]$, suppose that
\begin{align}
&\tag{$V_\alpha$} cR^\alpha\leq \lambda(B(x,R))\leq CR^\alpha \text{ for all } x\in G,\, R\geq1, \label{eq:Valpha}\\
&\tag{$G_\nu$} c\leq g(x,x)\leq C, \quad c\,d(x,y)^{-\nu}\leq g(x,y)\leq C\,d(x,y)^{-\nu} \text{ for
all } x\neq y\in G. \label{eq:Gnu}
\end{align}
We set $\beta\coloneqq\alpha-\nu$. 
The assumptions \eqref{eq:p_0}, \eqref{eq:Valpha} and \eqref{eq:Gnu} are very natural. Indeed for $d=d_G$ they were proved in landmark work of Grigor'yan-Telcs \cite{GT01} to be equivalent to two-sided sub-diffusive heat kernel bounds for $X$. In particular, these conditions imply that $E_x[T_{B(x,R)}]$ is of order $R^\beta$ for all $x\in G$ and $R\geq1$. Thus $\beta$ plays the role of the walk dimension. We refer to the introduction of \cite{drewitz_geometry_2024} for a number of examples of graphs that satisfy these conditions. Aside from regular lattices such as $\Z^d$, a case in point is the Cayley graph of the discrete Heisenberg group $\mathbb H_n(\Z)$ of dimension $2n+1$, with respect to its standard generators, which has polynomial growth of degree $2n+2$ and satisfies Gaussian heat kernel bounds \cite{hebisch_saloff-coste_1993}, so that \eqref{eq:Valpha} and \eqref{eq:Gnu} hold with $\alpha=2n+2$ and $\nu=2n$; for $n\geq2$ it fits the conditions of Theorem~\ref{thm:main} below, with $n=2$ corresponding exactly to the threshold $2\alpha=3\nu$ that will appear there. Note though that a notion of transitivity is not required, and that there are pertinent examples (e.g.~products of graphs) for which the natural choice of distance is $d \neq d_G$.

Let now
\begin{equation}\label{eq:theta}
\theta(x,R)\defeq\P\big(x\leftrightarrow B(x,R)^c\big),\qquad x\in G,\ R\geq1.
\end{equation}
Under the above assumptions, Proposition~\ref{prop:lb} has the following consequence (see Section~\ref{sec:prelim}).

\begin{corollary}\label{cor:lb}
Assume \eqref{eq:p_0}, \eqref{eq:Valpha} and \eqref{eq:Gnu}. Then for all $x\in G$ and $R\geq1$,
\begin{equation}\label{eq:cor_lb}
\theta(x,R)\geq c\max\big\{R^{-\nu/2},R^{-\beta}\big\}.
\end{equation}
\end{corollary}

 The exponents in \eqref{eq:cor_lb} come from the capacity and expected exit time of balls, see \eqref{eq:cap-ball} in Lemma~\ref{lem:ball} and the proof of Corollary~\ref{cor:lb} in Section~\ref{sec:prelim}. Both are sharp: in particular, the bound $E_x[T_{B}]\leq CR^\beta$, $B=B(x,R)$, used there is complemented by $E_x[T_{B}]\geq cR^\beta$, which follows from \eqref{eq:exit} as $1-E_x[e^{-\kappa T_B}]\leq\kappa E_x[T_B]$. 
 
 The bound \eqref{eq:cor_lb} is already insightful. The second (`volume') term on the right (see~Remark~\ref{rem:lb_general} as to the link with $\E|\cC_x\cap D|$) dominates if and only if $\beta \leq \nu/2$, or equivalently, $2\alpha\leq3\nu$.
 For $G=\Z^6$
($\alpha=6$, $\nu=4$), both terms equal $R^{-2}$. We will soon come back to this.

\medskip
Our next result is an upper bound. To state it, we introduce the quantity
\begin{equation}\label{eq:mod_triangle}
\tilde{\nabla}_{\kappa}\defeq\sup_{x\in G}\int_0^\infty\big(t\wedge \kappa^{-1}\big)^2P_x(X_t=x)\,dt,
\end{equation}
which is related to the (three-fold convolution) of the Green's function of the walk of $\mathcal{G}$ killed at rate $\kappa$; see Remark~\ref{rem:ub} below.

\begin{proposition}[Upper bound]\label{prop:ub}
Assume \eqref{eq:p_0} and $(\lambda_x\vee g(x,x))\leq C$ for all $x\in G$. Then for all $x\in G$, finite $D\ni x$, finite $B\subset
D^c$ and all $\kappa\in(0,1]$, 
\begin{equation}\label{eq:ub}
\P\big(x\leftrightarrow B\big)\leq C\,E_x\big[e^{-\kappa T_D}\,\P\big(X_{T_D}\leftrightarrow B\big)\big]+C\,\kappa\tilde{\nabla}_\kappa.
\end{equation}
\end{proposition}
 The assumptions on $\mathcal G$ are the same as in Proposition~\ref{prop:lb}. 
In practice, it is often sufficient to bound the connection probability in \eqref{eq:ub} by its maximum over $\partial D$, which yields
\begin{equation}\label{eq:ub-max}
\P\big(x\leftrightarrow B\big)\leq C\,E_x\big[e^{-\kappa T_D}\big]\max_{y\in\partial
D}\P\big(y\leftrightarrow B\big)+C\,\kappa\tilde{\nabla}_\kappa;
\end{equation}
this is the form in which Proposition~\ref{prop:ub} will be used to prove Theorem~\ref{thm:main}.

Proposition~\ref{prop:ub} is proved in Section~\ref{sec:ub}. It follows from two inputs. The first one compares $\P=\P_0$ and $\P_\kappa$, the law of the (centred) massive free field with mass $\kappa\lambda_x$ at vertex $x$ (see Section~\ref{sec:prelim} for the precise definition) on
$\sigma(\cC_x)$, at the price of an additive error term; it is proved in Section~\ref{sec:entropy} by
an entropy argument further discussed there. The second input is a Simon--Lieb-type inequality \cite{simon_correlation_1980,lieb_refinement_1980} for the killed field, proved in
Section~\ref{sec:SL}.

\begin{remark}[Triangle condition]\label{rem:ub} Under the assumptions \eqref{eq:p_0}, \eqref{eq:Valpha} and \eqref{eq:Gnu}, $\tilde{\nabla}_{\kappa}$ can be replaced by a more traditional triangle condition. 
For $\kappa\geq0$ and $x\in G$, let
\begin{equation}\label{eq:triangle}
\nabla_\kappa(x)\defeq\sum_{y,z\in
G}g_\kappa(x,y)\lambda_y\,g_\kappa(y,z)\lambda_z\,g_\kappa(z,x)=\tfrac12\int_0^\infty t^2e^{-\kappa
t}p_t(x,x)\,dt.
\end{equation}
Let $\nabla_\kappa\coloneqq\sup_{x\in G}\nabla_\kappa(x)$. One always has
$c\nabla_\kappa(x)\leq \tilde{\nabla}_{\kappa}\leq C\nabla_\kappa$ under
\eqref{eq:Valpha} and \eqref{eq:Gnu}; see Lemma~\ref{lem:kappa-tri} below.
The quantity $\tilde{\nabla}_\kappa$ is the one naturally produced by the entropy bound below, see Proposition~\ref{lem:entropy}.
\end{remark}

Let us say a few words on the choice of perturbation by a mass. The approach of \cite{dewan_muirhead_2023} (see also \cite{vEGPS}) for
Bernoulli percolation from which we drew inspiration (see also Remark~\ref{rem:sec3_end},\ref{rem:ber_compa} for a detailed comparison) perturbs the parameter $p$ itself, and the contraction step then requires
near-critical information at $p<p_c$, namely the sharp length estimate of
\cite{duminil-copin_alternative_2024,panis_schapira_2026}, which controls the quantity
$\varphi_p(S)$ of \cite{duminil-copin_tassion_2016}. The analogous perturbation for the free field
consists in raising the level, i.e.\ considering the clusters of $\{\tilde\varphi\geq\ell\}$ for
small $\ell>0$. The contraction step 
that an analogue of \eqref{eq:ub} with this perturbation would provide is then seen to involve no triangle condition at all, but rather to require a lower bound on the capacity of clusters realising the one-arm
event, in the spirit of the sharp length estimates of \cite{duminil-copin_alternative_2024}. This
looks like a promising avenue, which could also lead to meaningful off-critical results. We do not pursue it
here.

\medskip

As a consequence of {Proposition~\ref{prop:ub} and Corollary~\ref{cor:lb}}, we obtain the following result.

\begin{theorem}\label{thm:main}
Assume \eqref{eq:p_0}, \eqref{eq:Valpha} and \eqref{eq:Gnu}. If $2\alpha<3\nu$, then for all $x\in G$ and $R\geq1$,
\begin{equation}\label{eq:main}
cR^{-\beta}\leq\theta(x,R)\leq CR^{-\beta}.
\end{equation}
If $2\alpha=3\nu$, then $cR^{-\beta}\leq\theta(x,R)\leq CR^{-\beta}\log(1+R)$ for all $x\in G$ and
$R\geq1$.
\end{theorem}

The proof of Theorem~\ref{thm:main} appears in Section~\ref{sec:ub}. Let us now put this into context. Theorem~\ref{thm:main} complements the picture in lower dimensions. When $3\nu<2\alpha$, which
for $\cG=\Z^d$ means $d\in\{3,4,5\}$, the capacity term dominates in \eqref{eq:cor_lb} and the
one-arm probability is of order $R^{-\nu/2}$. This was proved on $\Z^3$ separately in \cite{dpr_low_2024,cai_ding_2024} (see also \cite{ding_wirth_2020} and \cite{dpr_intermediate_2023} for precursor results, with a $\sqrt{\log}$ and a $\log\log$ correction, respectively; and see \cite{dpr26fractal,CD_scaling} for recent related progress). Matching bounds for low-dimensional graphs with $\nu<1$ were previously derived in \cite{drewitz_critical_2023}. The argument of \cite{dpr_low_2024} is valid under \eqref{eq:p_0}, \eqref{eq:Valpha} and \eqref{eq:Gnu}, in combination with the lower bound of \cite{drewitz_critical_2023}: it yields bounds up to constants when $\nu<\alpha/2$ (i.e.\ up to dimension $4-\varepsilon$), and a polylogarithmic correction when $\nu=\alpha/2$. That of \cite{cai_ding_2024} on $\Z^d$ includes $d=4,5$, with a correction of order $e^{\sqrt{\log R}\,\log\log R}$ at $d=6$, and \cite{cai_one-arm_2025} treats $d>6$ by adapting the strategy of \cite{kozma_nachmias_2011}; high-dimensional ideas were already outlined in \cite{werner_clusters_2021}, to which we also refer for the conjectural picture across dimensions. The regime $2\alpha<3\nu$ was thus previously only treated for $\cG=\Z^d$, $d>6$. Theorem~\ref{thm:main} extends this to general graphs satisfying \eqref{eq:p_0}, \eqref{eq:Valpha} and \eqref{eq:Gnu}, with no transitivity assumption.

In the critical case
$2\alpha=3\nu$, Theorem~\ref{thm:main} sharpens the aforementioned bound of \cite{cai_ding_2024} on $\Z^6$ to a logarithmic correction. In
forthcoming work \cite{rodriguez_shi_zhang_prep}, we determine the precise order of $\theta(x,R)$
in this critical case and show, among other things, that for $G=\Z^6$ and all $R\geq2$,
\begin{equation}\label{eq:6d}
cR^{-2}\sqrt{\log R}\leq\theta(0,R)\leq CR^{-2}\sqrt{\log R}.
\end{equation}
Let us briefly comment on \eqref{eq:6d} (for a full heuristic we refer to \cite{rodriguez_shi_zhang_prep}). A simple but important observation is that at $2\alpha=3\nu$ , the system is simply {not} mean-field, it lies \emph{strictly below} the mean-field threshold. One can use the slab idea from \cite{RZ_cable}, used therein to interpolate precisely between dimensions two (the other critical dimension of the problem) and three, to see clearly where the transition to mean-field regime `sits'. Consider the slab $\Z^6\times(\Z/h\Z)$, for varying $h \geq 1$, which interpolates between $\Z^6$ and $\Z^7$: the heat kernel is that of $\Z^7$ for $t\leq h^2$ and of order $h^{-1}t^{-3}$ for $t \geq h^2$. Hence, with $\kappa=R^{-2}$ (i.e.~at spatial scale $R$), one readily finds that $\tilde\nabla_\kappa\approx 1 +h^{-1}\log R$, which remains bounded uniformly in $R$ precisely when $h \geq c\log R$, and the arguments of this article can be shown to extend to this regime, and to yield $\theta(x,R)\asymp R^{-2}$. In particular, $\Z^6$ lies strictly below the mean-field threshold (which occurs at $h \gg 1$). Crucially, this licenses the use of low-dimensional, capacity-governed heuristics (already implicit in \cite{dpr_low_2024,cai_ding_2024}) to describe the critical picture.

\medskip

The rest of this article is organised as follows. Section~\ref{sec:prelim} collects the facts
about the random walk, the free field and the loop soup that we need, and contains the proofs of
Proposition~\ref{prop:lb} and Corollary~\ref{cor:lb}. Section~\ref{sec:entropy} establishes the
entropy bound, Proposition~\ref{lem:entropy}, which compares the laws of $\cC_x$ under $\P$ and
$\P_\kappa$; the general inequality on which it rests, Lemma~\ref{lem:gen}, may be of independent
interest. Section~\ref{sec:SL} proves the Simon--Lieb-type inequality, Proposition~\ref{lem:SL}.
Section~\ref{sec:ub} combines these to prove Proposition~\ref{prop:ub} and Theorem~\ref{thm:main}.

\section{Preliminaries and lower bounds}\label{sec:prelim}

In this section, we first gather a few useful facts about the graphs $\cG$ under consideration and the behaviour of the walk $X$ on them, see Lemma~\ref{lem:ball}. We then compare the quantities $\nabla_\kappa$ and $\tilde\nabla_\kappa$, see Lemma~\ref{lem:kappa-tri}. Finally, we collect a few facts about the free field, including their dependence on a mass $\kappa \geq 0$, and prove the lower bounds Proposition~\ref{prop:lb} and Corollary~\ref{cor:lb}. Below, for two positive functions $f$ and $g$ with the same domain, we write $f\asymp g$ if $c\leq f/g\leq C$.

Recall the setup from the beginning of Section~\ref{sec:intro}, by which $\mathcal G=(G,\lambda)$ is a transient weighted graph, connected and locally finite, and $X$ is the associated random walk with unit jump rate. Any additional assumption on $\mathcal{G}$ will always be made explicit.

\begin{lemma}\label{lem:ball}
Assume \eqref{eq:p_0}, \eqref{eq:Valpha} and \eqref{eq:Gnu}. 
Then for all $x\in G$ and $R\geq1$, writing $B=B(x,R)$,
\begin{align}
&\capacity(B)\asymp R^{\nu}, \label{eq:cap-ball}\\
&\textstyle\sum_{y\in B}g(x,y)\asymp R^\beta, \label{eq:gsum}\\
&E_x\big[e^{-\kappa T_B}\big]\leq C\exp\big\{-c(\kappa R^\beta)^{1/\beta}\big\},\quad\text{$\kappa\in(0,1]$}, \label{eq:exit}\\
&p_t(x,x)\asymp t^{-\alpha/\beta},\quad\text{$t\geq1$}. \label{eq:hk}
\end{align}
\end{lemma}

\begin{proof}
The bound \eqref{eq:cap-ball} is \cite[(3.11)]{drewitz_geometry_2024}. For \eqref{eq:gsum}, by \eqref{eq:Valpha}, \eqref{eq:Gnu} and $\lambda_y\asymp1$ (see \cite[(2.10)]{drewitz_geometry_2024}), summing over the dyadic shells $B(x,2^k)\setminus B(x,2^{k-1})$ gives $\sum_{y\in B}g(x,y)\leq C\sum_{0\leq k\leq \log_2 R+1}2^{-k\nu}2^{k\alpha}\leq CR^\beta$, while $\sum_{y\in B}g(x,y)\geq cR^{-\nu}|B|\geq cR^{-\nu}\lambda(B)\geq cR^\beta$. The remaining bounds are standard, see \cite[Chap.~4--5]{Barlow_2017}, but are stated there for $d=d_G$; we reproduce the arguments to highlight that they do not depend on this choice. Let $Z$ be the discrete skeleton of $X$ and $N_t$ the number of jumps of $X$ up to time $t$, so that $X_t=Z_{N_t}$ with $N$ a unit-rate Poisson process independent of $Z$. By \cite[Prop.~3.3]{drewitz_geometry_2024}, $P_x(Z_n=y)\leq C\lambda_y n^{-\alpha/\beta}$ and $P_x(T^Z_B\leq n)\leq C\exp\{-(cR^\beta/n)^{1/(\beta-1)}\}$ for $n\geq1$, where $T^Z_B$ is the exit time of $B$ for $Z$. Moreover, $P(|N_t-t|\geq t/2)\leq 2e^{-ct}$, see e.g.\ \cite[Lemma~5.13]{Barlow_2017}. Since $T_B\leq t$ if and only if $N_t\geq T^Z_B$, it follows that for $t\geq R$,
\begin{equation}\label{eq:exit-tail}
\begin{split}
P_x\big(T_B\leq t\big)&\leq P(N_t\geq2t)+P_x\big(T^Z_B\leq\lfloor2t\rfloor\big)\\
&\leq 2e^{-ct}+C\exp\big\{-(cR^\beta/2t)^{1/(\beta-1)}\big\}\leq C\exp\big\{-(cR^\beta/t)^{1/(\beta-1)}\big\},
\end{split}
\end{equation}
where we used that $(R^\beta/t)^{1/(\beta-1)}\leq t$ for $t\geq R$, cf.\ \cite[Lemma~5.21]{Barlow_2017}. For \eqref{eq:exit}, write $E_x\big[e^{-\kappa T_B}\big]\leq P_x(T_B\leq t)+e^{-\kappa t}$
and 
take $t=R\kappa^{-(\beta-1)/\beta}\ (\geq R)$ in \eqref{eq:exit-tail} to get the result.

We now turn to \eqref{eq:hk}. For the upper bound, as $\lambda_x\asymp1$, $p_t(x,x)\leq\lambda_x^{-1}P(N_t<t/2)+C(t/2)^{-\alpha/\beta}\leq Ct^{-\alpha/\beta}$ for $t\geq1$. For the lower bound we follow \cite[Lemma~4.35]{Barlow_2017}, see also \cite{CoulhonGrigoryan97}: by symmetry of $p_t$, the semigroup property, the Cauchy--Schwarz inequality and \eqref{eq:Valpha},
\[
p_{2t}(x,x)=\sum_{y\in G}\lambda_y\,p_t(x,y)^2\geq\lambda(B)^{-1}P_x(X_t\in B)^2\geq cR^{-\alpha}\big(1-P_x(T_B\leq t)\big)^2.
\]
Choose $R=(\CL{C:hk}t)^{1/\beta}$, where $\CR{C:hk}\geq1$ is large enough that the right side of \eqref{eq:exit-tail} is at most $\frac12$ when $R^\beta/t=\CR{C:hk}$. Then \eqref{eq:exit-tail} yields $P_x(T_B\leq t)\leq\frac12$ for $t\geq t_0\coloneqq\CR{C:hk}^{1/(\beta-1)}$ (which ensures $t\geq R$), whence $p_{2t}(x,x)\geq ct^{-\alpha/\beta}$. Meanwhile for $t\in[1,2t_0]$, one bounds $p_t(x,x)\geq\lambda_x^{-1}e^{-t}\geq c$.
\end{proof}

We now turn to the relationship between $\tilde \nabla_{\kappa}$ and $\nabla_{\kappa}$ introduced in \eqref{eq:mod_triangle} and below \eqref{eq:triangle}, respectively. Let $\epsilon\coloneqq3\nu-2\alpha\ (=\alpha-3\beta)$ and the non-decreasing function
\begin{equation}\label{eq:def-b}
b(R)\defeq 1+\int_1^Rr^{-1-\epsilon}\,dr,\quad R\geq1,
\end{equation}
so that $b(R)\asymp1$ if $\epsilon>0$, $b(R)\asymp\log(1+R)$ if $\epsilon=0$ and $b(R)\asymp R^{-\epsilon}$ if $\epsilon<0$, for all $R\geq1$.

\begin{lemma}\label{lem:kappa-tri}
Assume \eqref{eq:p_0}, \eqref{eq:Valpha} and \eqref{eq:Gnu}. For all $\kappa\in(0,1]$, one has
\begin{equation}\label{eq:kappa-tri}
c\,b(\kappa^{-1/\beta})\leq\nabla_\kappa\leq C\tilde\nabla_\kappa\leq C^\prime b(\kappa^{-1/\beta}).
\end{equation}
\end{lemma}

\begin{proof}
Write $\tilde\nabla_\kappa(x)\coloneqq \lambda_x \int_0^\infty(t\wedge\kappa^{-1})^2p_t(x,x)\,dt$, so that $\tilde\nabla_\kappa=\sup_x\tilde\nabla_\kappa(x)$, and recall that $\lambda_x\asymp1$, see \cite[(2.10)]{drewitz_geometry_2024}. Since $s^2e^{-s}\leq1$ for $s\geq0$, we have $t^2e^{-\kappa t}\leq(t\wedge\kappa^{-1})^2$, whence $\nabla_\kappa(x)\leq C\tilde\nabla_\kappa(x)$. Next, by the upper bound in \eqref{eq:hk}, $\alpha>\beta$ and the substitution $t=r^\beta$,
\[
\tilde\nabla_\kappa(x)\leq C\int_0^1t^2\,dt+C\int_1^{\kappa^{-1}}t^{2-\alpha/\beta}\,dt+C\kappa^{-2}\int_{\kappa^{-1}}^\infty t^{-\alpha/\beta}\,dt\leq C\big(b(\kappa^{-1/\beta})+\kappa^{\alpha/\beta-3}\big)\leq Cb(\kappa^{-1/\beta}),
\]
as $\kappa^{\alpha/\beta-3}=R^{-\epsilon}$ for $R=\kappa^{-1/\beta}$. Finally, as $e^{-\kappa t}\geq e^{-1}$ for $t\leq\kappa^{-1}$, $P_x(X_t=x)\geq e^{-t}$ and by the lower bound in \eqref{eq:hk},
\[
\nabla_\kappa(x)\geq c\int_0^1t^2\,dt+c\int_1^{\kappa^{-1}}t^{2-\alpha/\beta}\,dt\geq c\,b(\kappa^{-1/\beta}).
\]
Taking suprema over $x$ yields \eqref{eq:kappa-tri}.
\end{proof}

We now turn to matters relating to the Gaussian free field, and allow for reasons to become clear later for the possibility of a small mass $\kappa\in[0,1]$. In order to accommodate this within our setup, for $\kappa\geq0$, let $\cG_\kappa\coloneqq(G,\lambda,\kappa\lambda)$ denote the graph $\cG$ endowed with the killing measure $\kappa_x\coloneqq\kappa\lambda_x$, $x\in G$, see \cite[Section~2]{drewitz_cluster_2022}; the choice of $\kappa_x$ proportional to $\lambda_x$ corresponds to killing $X$ at constant rate $\kappa$, so that the Green function of $\cG_\kappa$ is $g_\kappa$ from \eqref{eq:green}. Its metric graph $\tilde\cG_\kappa$ is obtained from $\tilde\cG$ by attaching to each $x\in G$ an additional half-open cable of length $\rho_x\coloneqq(2\kappa_x)^{-1}$. Let $\P_\kappa$ denote the law of the free field on $\tilde\cG_\kappa$, so $\P_0=\P$. With $\rho_x=\infty$ when $\kappa=0$, $\tilde\cG_0$ is $\tilde\cG$ with a half-open infinite cable attached at each vertex, which does not affect the law of the free field on $\tilde\cG$, cf.\ \cite[Section~2]{drewitz_cluster_2022}. In this way, we view $\tilde\cG_\kappa$  as a subset of $\tilde\cG_0$ for all $\kappa>0$, the diffusion on $\tilde\cG_\kappa$ being the diffusion on $\tilde\cG_0$ killed upon reaching the open end of the cable of length $\rho_x$ at some $x\in G$. This is particularly handy for monotonicity statements involving the loop soup, see \eqref{eq:mono} below. We write $\capacity_\kappa$ for the capacity associated to $\cG_\kappa$ (so $\capacity_0=\capacity$), see \cite[Section~2]{drewitz_cluster_2022}; it extends to compact subsets of $\tilde\cG_\kappa$, it is monotone, and for finite $K\subset G$ the union $\tilde K$ of $K$ and of all cables with both endpoints in $K$ satisfies $\capacity_\kappa(\tilde K)=\capacity_\kappa(K)$. If $\sup_xg_\kappa(x,x)<\infty$, then for all $x\in G$ and $t>g_\kappa(x,x)^{-1}$,
\begin{equation}\label{eq:caplaw}
\P_\kappa\big(\capacity_\kappa(\cC_x)\geq t\big)=\frac1\pi\arctan\big((g_\kappa(x,x)t-1)^{-1/2}\big),
\end{equation}
see \cite[Theorem~3.7]{drewitz_cluster_2022}. Moreover, for all $x,y\in G$, see
\cite{lupu_loop_2016},
\begin{equation}\label{eq:twopoint}
\P_\kappa(y\in\cC_x)=\frac1\pi\arcsin\Big(\frac{g_\kappa(x,y)}{\sqrt{g_\kappa(x,x)g_\kappa(y,y)}}\Big).
\end{equation}
 We will in fact only ever use the lower bound on the tail of $\capacity(\cC_x)$ implied by \eqref{eq:caplaw} (for $\kappa=0$), and the upper bound on $\P_\kappa(y\in\cC_x)$ implied by \eqref{eq:twopoint}. 

Without any assumption on $\cG$, as $g_\kappa(x,x)\geq((1+\kappa)\lambda_x)^{-1}$ (cf.\ below \eqref{eq:lb}), \eqref{eq:twopoint} and the inequality $\arcsin u\leq\frac\pi2u$ for $u\in[0,1]$ yield $\P_\kappa(y\in\cC_x)\leq\tfrac12g_\kappa(x,y)(g_\kappa(x,x)g_\kappa(y,y))^{-1/2}\leq g_\kappa(x,y)(\lambda_x\lambda_y)^{1/2}$ for $\kappa\in[0,1]$. As $\sum_{y\in G}\lambda_yg_\kappa(x,y)=\int_0^\infty e^{-\kappa t}\,dt=\kappa^{-1}$, it follows that 
\begin{equation}\label{eq:vol_bound_general}
\E_\kappa\big[|\cC_x\cap G|\big]\leq\sum_{y\in G}\lambda_yg_\kappa(x,y)\Big(\frac{\lambda_x}{\lambda_y}\Big)^{1/2}\leq\kappa^{-1}\sup_{y\in G}\Big(\frac{\lambda_x}{\lambda_y}\Big)^{1/2},
\end{equation}
for all $x\in G$ and $\kappa\in(0,1]$.
This finishes our preparation.
We now prove Proposition~\ref{prop:lb} and Corollary~\ref{cor:lb}, which only concern the case $\kappa=0$.

\begin{proof}[Proof of Proposition~\ref{prop:lb}]
 Let $x\in G$ and $D\ni x$ finite. 
By \eqref{eq:caplaw}, the inequality $\arctan u\geq\frac\pi4u$ valid for $u\in[0,1]$ and the bound $g(x,x)\leq C$, we get
\begin{equation}\label{eq:cap_tail}
\P\big(\capacity(\cC_x)\geq t\big)\geq\tfrac14\big(g(x,x)t\big)^{-1/2}\geq ct^{-1/2},\quad t\geq2g(x,x)^{-1}.
\end{equation}

\textit{Capacity term.} On $\{x\nlr D^c\}$, all vertices of $\cC_x$ lie in $D$, whence $\cC_x$ is contained in the union of $\overline D\coloneqq D\cup\partial D$ and of all cables with both endpoints in $\overline D$, and so by monotonicity,
\begin{equation}\label{eq:cap_D}
\capacity(\cC_x)\leq\capacity(\overline D)\leq c_0^{-2}\capacity(D)\quad\text{on } \{x\nlr D^c\}.
\end{equation}
To get the second inequality, $h(y)\coloneqq P_y(H_D<\infty)\geq c_0$ for
$y\in\partial D$ by \eqref{eq:p_0}, so $(h/c_0)\wedge1$ is admissible in the Dirichlet principle for
$\capacity(\overline D)$, see \cite[(1.62)]{sznitman_topics_2012}, and its energy is at most
$c_0^{-2}$ times that of $h$, which equals $\capacity(D)$. Let $t_1\coloneqq2c_0^{-2}\capacity(D)$, so that $t_1\geq2g(x,x)^{-1}$ as $\capacity(D)\geq\capacity(\{x\})=g(x,x)^{-1}$. By \eqref{eq:cap_D} and \eqref{eq:cap_tail}, $\P(x\leftrightarrow D^c)\geq\P(\capacity(\cC_x)\geq t_1)\geq c\,\capacity(D)^{-1/2}$.

\medskip

\textit{Volume term.} Let $\P^D$ denote the law of the free field on $\tilde\cG$ with Dirichlet boundary condition on $D^c$, i.e.~of the centred Gaussian field on $\tilde\cG$ vanishing on $D^c$ whose covariance restricted to $D$ is the killed Green's function $
g^D(y,z)\coloneqq\lambda_z^{-1}E_y\big[\int_0^{T_D}1\{X_t=z\}\,dt\big]$,
and let $\cC_x^D$ denote the cluster of $x$ in $\{\tilde\varphi>0\}$ under $\P^D$. Using the coupling of loop soup clusters and the free field of \cite[Proposition~2.1]{lupu_loop_2016}, and restricting the loop soup to the loops avoiding $D^c$, one can couple $\cC_x$ under $\P$ and $\cC_x^D$ under $\P^D$ in such a way that $\cC_x^D\subset\cC_x$, with equality on $\{x\nlr D^c\}$. Consequently, for all $t>0$,
\begin{equation}\label{eq:loc_D}
\P\big(x\leftrightarrow D^c\big)\geq\P\big(\capacity(\cC_x)\geq t\big)-\P^D\big(\capacity(\cC_x^D)\geq t\big).
\end{equation}
 Let $K\coloneqq\overline{\cC_x^D\cap G}$, where for $D \subset G$ we write $\overline{D}= D\cup \partial D$.
As $\cC_x^D\subset\tilde K$ (see above \eqref{eq:caplaw} for notation) and \eqref{eq:p_0} bounds the degree by $c_0^{-1}$, subadditivity of the capacity and $g(y,y)\geq c$ give
\begin{equation}\label{eq:cap_vol}
\capacity(\cC_x^D)\leq\sum_{y\in K}g(y,y)^{-1}\leq C|\cC_x^D\cap G|.
\end{equation}
By a calculation similar to the one leading to \eqref{eq:vol_bound_general}, using the analogue of \eqref{eq:twopoint} under $\P^D$, see \cite{lupu_loop_2016,drewitz_cluster_2022}, as well as the lower bounds $g^D(y,y)\geq\lambda_y^{-1}\geq c$ and $\lambda_y\geq c$, we get
\begin{equation}\label{eq:vol_D}
\E^D\big[|\cC_x^D\cap G|\big]\leq\sum_{y\in D}\frac{g^D(x,y)}{2\sqrt{g^D(x,x)g^D(y,y)}}\leq C\sum_{y\in D}\lambda_yg^D(x,y)=C\,E_x[T_D].
\end{equation}
By \eqref{eq:loc_D}, \eqref{eq:cap_tail}, \eqref{eq:cap_vol}, Markov's inequality and \eqref{eq:vol_D}, for $t\geq C( \geq 2g(x,x)^{-1})$,
\[
\P\big(x\leftrightarrow D^c\big)\geq ct^{-1/2}-C't^{-1}E_x[T_D].
\]
Choosing $t=(A\,E_x[T_D])^2$, which is admissible for $A\geq C$ since $E_x[T_D]\geq1$, 
yields $\P(x\leftrightarrow D^c)\geq(cA^{-1}-C'A^{-2})E_x[T_D]^{-1}\geq c'E_x[T_D]^{-1}$ upon taking $A$ large. Overall, \eqref{eq:lb} follows.\black
\end{proof}

\begin{proof}[Proof of Corollary~\ref{cor:lb}]
 Apply \eqref{eq:lb} with $D=B(x,R)$ (Proposition~\ref{prop:lb} applies as $\lambda_x\asymp1$, see \cite[(2.10)]{drewitz_geometry_2024}), and use \eqref{eq:cap-ball} together with $E_x[T_B]\leq\sum_{y\in B}\lambda_yg(x,y)\leq CR^\beta$, valid by \eqref{eq:gsum}.\black
\end{proof}

\begin{remark}[Weaker assumptions]\label{rem:lb_general}
Under \eqref{eq:p_0} and $c\leq g(x,x)\leq C$ for all $x\in G$ only (which allows $\lambda$ to be unbounded from above), the above argument applied directly to $\cC_x$ instead of $\cC_x^D$ yields, for all $x\in G$ and finite $D\ni x$,
\begin{equation}\label{eq:lb_green}
\P\big(x\leftrightarrow D^c\big)\geq c\max\Big\{\capacity(D)^{-1/2},\Big(\sum_{y\in
D}g(x,y)\Big)^{-1}\Big\}.
\end{equation}
Indeed, on $\{x\nlr D^c\}$, one has $\capacity(\cC_x)\leq c_0^{-2}\capacity(D)$ and $\capacity(\cC_x)\leq C|\cC_x\cap D|$ as in \eqref{eq:cap_D} and \eqref{eq:cap_vol}. The first bound gives the capacity term as above. For the second term, one applies Markov's inequality directly to the volume, using the bound $\E[|\cC_x\cap D|]\leq C\sum_{y\in D}g(x,y)$, which follows using \eqref{eq:twopoint} and the assumptions on $\cG$. Under the assumptions of Proposition~\ref{prop:lb}, \eqref{eq:lb} implies \eqref{eq:lb_green}, as $E_x[T_D]=\sum_{y\in D}\lambda_yg^D(x,y)\leq C\sum_{y\in D}g(x,y)$. The gain can potentially be substantial, e.g.\ when $x$ is close to $\partial D$.
\black
\end{remark}

\section{Entropy bound}\label{sec:entropy}

In this section, we prove an entropy bound, Proposition~\ref{lem:entropy}, which compares the law of the (massless) cluster $\cC_x$ under $\P=\P_0$ with its massive counterpart (under $\P_\kappa$). It is inspired by a method of \cite{dewan_muirhead_2023}; see also \cite{hutchcroft_2022_on, vEGPS} in the context of Bernoulli percolation. Here we perturb by introducing a mass rather than changing the level (or, equivalently, the density $p$). 

The proof has two steps. First, Lemma~\ref{lem:gen} gives a general inequality that bounds $\P(A)$ by $(1+\varepsilon)\P_\kappa(A)$ plus a relative entropy term; this is effectively a Gibbs variational principle/Donsker-Varadhan  type argument (see Remark~\ref{rem:pinsker},\ref{rem:DV} below). Second, we bound this relative entropy in the loop soup picture. The loop soup is the natural language for measuring the effect of $\kappa$, which re-weighs loops exponentially in their duration. The relative entropy is computed for general Poisson processes (Lemma~\ref{lem:poisson}), with the flexibility to perform an (adapted) exploration argument. This allows to progressively reveal the loops that make up the cluster $\cC_x$. Via a union bound over the vertices visited by these loops, the entropy bound reduces to `one-point' computations for loops through a single vertex $y\in G$, see Lemma~\ref{lem:one_point}. 

\black
Recall $\tilde{\nabla}_{\kappa}$ from \eqref{eq:mod_triangle}, the graph $\cG_\kappa$, its metric graph $\tilde\cG_\kappa$ and the massive free field $\P_\kappa$ from Section~\ref{sec:prelim}. The main result of this section is the following:

\begin{proposition}[Entropy bound]\label{lem:entropy}
Assume  \eqref{eq:p_0}. \black Then for all $x\in G$, $\kappa\in(0,1]$ and
$A\in\sigma(\cC_x)$,
\begin{equation}\label{eq:entropy}
\P(A)\leq2\,\P_\kappa(A)+C{\kappa^2\chi_\kappa(x)}\,\tilde{\nabla}_{\kappa}, \quad \text{where }\chi_\kappa(x)\defeq\E_\kappa\big[|\cC_x\cap G|\big].
\end{equation}
\end{proposition}

The factor of $2$ in \eqref{eq:entropy} is arbitrary: the proof yields $\P(A)\leq(1+\varepsilon)\,\P_\kappa(A)+C\varepsilon^{-1}\kappa^2\chi_\kappa(x)\black\tilde{\nabla}_{\kappa}$ for all $\varepsilon\in(0,1]$.  Under the assumptions of Proposition~\ref{prop:ub}, the susceptibility will be controlled via \eqref{eq:vol_bound_general}, which yields $\chi_\kappa(x)\leq C\kappa^{-1}$ and leads to the error term in \eqref{eq:ub}, see Section~\ref{sec:ub}.\black

Note that \eqref{eq:entropy} is trivial for $\kappa=0$. In the other direction, one has the following. Call an event $A\in\sigma(\cC_x)$ \emph{increasing} if $A=\{\cC_x\in\mathcal{K}\}$ for some collection $\mathcal{K}$ of subsets of $\tilde\cG$ closed under taking supersets, i.e.\ $K\in\mathcal{K}$ and $K\subset K'$ imply $K'\in\mathcal{K}$ (we will only use this in the specific case $A= \{x\leftrightarrow B\}$). Then for all $x\in G$, $0\leq\kappa\leq\kappa'$ and increasing $A\in\sigma(\cC_x)$,
\begin{equation}\label{eq:mono}
\P_{\kappa'}(A)\leq\P_\kappa(A) \, (\leq \P(A)).
\end{equation}
Indeed, viewing $\tilde\cG_{\kappa'}\subset\tilde\cG_\kappa$ as in Section~\ref{sec:prelim}, the loop soup on $\tilde\cG_{\kappa'}$ is obtained from that on $\tilde\cG_\kappa$ by keeping only the loops contained in $\tilde\cG_{\kappa'}$, and \eqref{eq:mono} follows by the isomorphism of \cite{lupu_loop_2016} (see also \eqref{eq:loop_gff} below). \black

\medskip

Towards \eqref{eq:entropy},
we start by recording a general inequality relating the probabilities of an event under two measures to
their relative entropy, which is of independent interest. For probability measures $P,Q$ on a common measurable space, let $H(Q|P)\coloneqq E^Q[\log\frac{dQ}{dP}]$
if $Q\ll P$ and $H(Q|P)\coloneqq\infty$ otherwise, and introduce (with the convention $0\cdot (-\infty) =0$, so $\ophi(0)=1$)
\begin{equation}\label{eq:phi}
\ophi(u)\defeq u\log u-u+1,\quad u\geq0
\end{equation}
Then $\ophi$ is convex, non-negative, decreasing on $[0,1]$, $\ophi(1)=0$, and since
$E^P[\frac{dQ}{dP}]=1$,
\begin{equation}\label{eq:Hphi}
H(Q|P)=E^P\Big[\ophi\Big(\frac{dQ}{dP}\Big)\Big]\qquad\text{if }Q\ll P.
\end{equation}

\begin{lemma}\label{lem:gen}
If $Q\ll P$ and $P(A)>0$, then $P(A)\,\ophi\big(Q(A)/P(A)\big)\leq H(Q|P)$. In particular, for any event $A$ and all $\varepsilon\in(0,1]$,
\begin{equation}\label{eq:gen}
P(A)\leq(1+\varepsilon)\,Q(A)+2\varepsilon^{-1}H(Q|P).
\end{equation}
\end{lemma}
\begin{proof}
Combining \eqref{eq:Hphi}, the fact that $\ophi\geq0$ and Jensen's inequality, we get that
\[
H(Q|P)\geq E^P\Big[\ophi\Big(\frac{dQ}{dP}\Big)1_A\Big]=P(A)\,E^P\Big[\ophi\Big(\frac{dQ}{dP}\Big)\,\Big|\,A\Big]
\geq P(A)\,\ophi\Big(E^P\Big[\frac{dQ}{dP}\,\Big|\,A\Big]\Big),
\]
and $E^P[\frac{dQ}{dP}\,|\,A]=Q(A)/P(A)$.
For \eqref{eq:gen}, note that the Legendre transform $ \ophi^* (s)\coloneqq\sup_{u\geq0}\{su-\ophi(u)\} = e^s-1$ for all $s\in\R$ (the supremum being attained at $u=e^s$). Hence, using that $su\leq \ophi(u)+ \ophi^*(s) =\ophi(u)+e^s-1$, for all $s\in\R$ and $u\geq0$, with $s=-\varepsilon$ and $u=Q(A)/P(A)$, multiplying by $P(A)$ and using the first part, we obtain
\[
-\varepsilon\,Q(A)\leq P(A)\,\ophi\big(Q(A)/P(A)\big)-(1-e^{-\varepsilon})P(A)\leq H(Q|P)-(1-e^{-\varepsilon})P(A).
\]
 Dividing by $1-e^{-\varepsilon}$, using that $\frac{\varepsilon}{1-e^{-\varepsilon}}\leq1+\varepsilon$ and $\frac{1}{1-e^{-\varepsilon}}\leq\frac2\varepsilon$ for $\varepsilon\in(0,1]$, \eqref{eq:gen} follows.
\end{proof}

\begin{remark}\label{rem:pinsker}
\begin{enumerate}[label=\arabic*),itemsep=4pt]
\item\label{rem:DV} Readers familiar with the theory of large deviations will recognise in the bound $P(A)\,\ophi\big(Q(A)/P(A)\big)\leq H(Q|P)$ a consequence of the Donsker--Varadhan variational formula, see e.g.\ \cite[Lemma~6.2.13]{dembo_zeitouni_1998}. In fact, the latter leads to a better estimate with the left-hand side replaced by $\sum_{E\in \{A,A^c\}} P(E)\,\ophi\big(Q(E)/P(E)\big)$, which also follows directly from the proof of Lemma~\ref{lem:gen} by applying Jensen on $A^c$ as well.

\item Since $\ophi(1)=\ophi'(1)=0$ and $\ophi''(s)=s^{-1}$, Taylor's formula gives $\ophi(u)=\int_1^u\frac{u-s}{s}\,ds\geq\frac{(1-u)^2}{2\max\{1,u\}}$ for all $u\geq0$. Inserted into $P(A)\,\ophi\big(Q(A)/P(A)\big)\leq H(Q|P)$, this immediately yields the generalised Pinsker-type inequality of \cite[Lemma~2.3]{dewan_muirhead_2023}, see also \cite[(3.2)]{hutchcroft_2022_on}. 

One can also envision (in the spirit of \cite{dewan_muirhead_2023,hutchcroft_2022_on}) to prove a version of \eqref{eq:gen} by replacing the use of Young's inequality by distinguishing whether $Q(A)\geq P(A)/(1+\varepsilon)$ or not (and combining with the first inequality in Lemma~\ref{lem:gen}. Owing to $\ophi\big(\tfrac1{1+\varepsilon}\big)\geq c{\varepsilon^2}$, this leads to a slightly worse bound with $C\varepsilon^{-2}H(Q|P)$ in place of $2\varepsilon^{-1}H(Q|P)$.
\end{enumerate}
\end{remark}

In view of \eqref{eq:gen}, applied with $\varepsilon=1$, $P=\P$ and $Q=\P_\kappa$, the proof of Proposition~\ref{lem:entropy} boils down to bounding the relative entropy of the law of $\cC_x$ under $\P_\kappa$ relative to $\P$. To do this, we will regard the sign clusters as generated by Poissonian loop soups, see \eqref{eq:loop_gff} below. Hence, we focus on computing relative entropies for Poisson processes. Much of this is completely general.

Let $(S,\mathcal{S})$ be a measurable space. For a $\sigma$-finite measure $\mu$ on $S$, let $P_\mu$
denote the law of a Poisson point process $\omega$ on $S$ with intensity measure $\mu$, viewed as a
random point measure on $S$; we write $s\in\omega$ for the atoms of $\omega$, counted with
multiplicity. For $S'\in\mathcal{S}$, let $\omega^{S'}\coloneqq\omega|_{S'}$; under $P_\mu$,
$\omega^{S'}$ has law $P_{\mu|_{S'}}$ and is independent of $\omega^{S\setminus S'}$.

\begin{lemma}\label{lem:poisson}
Let $\mu_0$ be $\sigma$-finite, $\rho:S\to[0,1]$ measurable, and $\mu_1\coloneqq\rho\,\mu_0$.
\begin{itemize}
\item 
 If  $\int(1-\rho)\,d\mu_0<\infty$, then $P_{\mu_1}\ll P_{\mu_0}$, and
\begin{equation}\label{eq:poisson_entropy}
H\big(P_{\mu_1}\,\big|\,P_{\mu_0}\big)=\int\ophi(\rho)\,d\mu_0\ \Big(\leq\int(1-\rho)\,d\mu_0\Big).
\end{equation}
\item Let $S_0, S_1,\dots$ be random, pairwise disjoint elements of
$\mathcal{S}$ such that
$S_0$ is deterministic. Let $T_k=(\omega^{S_0},\dots,\omega^{S_{k}})$ and $\mathcal{F}_k=\sigma(T_k)$. Suppose
$S_n$ is $\mathcal{F}_{n-1}$-measurable
for all $n\geq1$, and $\int_{S_n}(1-\rho)\,d\mu_0<\infty$ $P_{\mu_1}$-a.s.\ for all $n\geq0$. Then for
$U_n\coloneqq\bigcup_{i=0}^nS_i$ and $U_\infty\coloneqq\bigcup_{i=0}^\infty S_i$,
\begin{equation}\label{eq:explore}
H\Big(\Law_{\mu_1}(\omega^{U_\infty})\,\big|\,\Law_{\mu_0}(\omega^{U_\infty})\Big)=E_{\mu_1}\Big[\int_{U_\infty}\ophi(\rho)\,d\mu_0\Big].
\end{equation}
\end{itemize}
\end{lemma}
\begin{proof}
We first show that
\begin{equation}\label{eq:poisson_density}
\frac{dP_{\mu_1}}{dP_{\mu_0}}(\omega)=f_S(\omega),\quad \text{where }f_S(\omega)\defeq\exp\Big\{\int_S(1-\rho)\,d\mu_0\Big\}\prod_{s\in\omega}\rho(s).
\end{equation}
 Recall that $E_\mu[e^{-\langle\omega,v\rangle}]=\exp\{-\int(1-e^{-v})\,d\mu\}$ for all measurable $v:S\to[0,\infty]$, see \cite[Thm.~3.9]{last_penrose_2017} (for $[0,\infty]$-valued $v$ by monotone convergence). Applying this with $v=u-\log\rho$, we obtain for all measurable $u:S\to[0,\infty)$, as $\int(1-\rho)\,d\mu_0<\infty$,
\begin{align*}
E_{\mu_0}\big[f_S\,e^{-\langle\omega,u\rangle}\big]&=\exp\Big\{\int(1-\rho)\,d\mu_0\Big\}E_{\mu_0}\big[e^{-\langle\omega,u-\log\rho\rangle}\big]\\
&=\exp\Big\{\int(1-\rho)\,d\mu_0-\int(1-\rho e^{-u})\,d\mu_0\Big\}=E_{\mu_1}\big[e^{-\langle\omega,u\rangle}\big].
\end{align*}
For $u=0$, this shows that $f_S\,dP_{\mu_0}$ is a probability measure, and \eqref{eq:poisson_density} follows since the Laplace functional determines the law of a point process, see \cite[Prop.~2.10]{last_penrose_2017}.
Now, $P_{\mu_1}$-a.s.\ $\rho(s)>0$ for all $s\in\omega$, and $\rho\log\rho\geq\rho-1$ is $\mu_0$-integrable. Hence, by Campbell's formula \cite[Proposition~2.7]{last_penrose_2017},
\begin{align*}
H\big(P_{\mu_1}\,\big|\,P_{\mu_0}\big)&=E_{\mu_1}[\log f_S]=\int(1-\rho)\,d\mu_0+E_{\mu_1}\Big[\sum_{s\in\omega}\log\rho(s)\Big]\\
&=\int(1-\rho+\rho\log\rho)\,d\mu_0=\int\ophi(\rho)\,d\mu_0,
\end{align*}
and the bound in parentheses in \eqref{eq:poisson_entropy} follows from $\rho\log\rho\leq0$. This proves \eqref{eq:poisson_entropy}.

It remains to show \eqref{eq:explore}. Since $S_k$ is $\mathcal F_{k-1}$-measurable and disjoint from $U_{k-1}$, conditionally on $\mathcal F_{k-1}$, $\omega^{S_k}$ is a Poisson  process under $P_{\mu_i}$ with intensity $\mu_i|_{S_k}$, $i\in\{0,1\}$, and $\int_{S_k}(1-\rho)\,d\mu_0<\infty$ $P_{\mu_1}$-a.s. Hence, by the chain rule for the relative entropy, see \cite[Section~5.3]{gray_1990_entro},\black\ and \eqref{eq:poisson_entropy} applied conditionally on $\mathcal F_{k-1}$, induction on $k$ gives, for all $k\geq0$,
\begin{equation}\label{eq:explore_after_iterate}
H\big(\Law_{\mu_1}(T_k)\,\big|\,\Law_{\mu_0}(T_k)\big)=\sum_{j=0}^k E_{\mu_1}\Big[\int_{S_j}\ophi(\rho)\,d\mu_0\Big]= E_{\mu_1}\Big[\int_{U_k}\ophi(\rho)\,d\mu_0\Big].
\end{equation}
As $k\to\infty$, the right side converges to that of \eqref{eq:explore} by monotone convergence, as $\ophi\geq0$, and the left-hand side converges to that of \eqref{eq:explore} by \cite[Corollary~5.2.3]{gray_1990_entro}, as $\sigma(T_k)\uparrow\sigma(\omega^{U_\infty})$ (this can also be deduced from the Donsker--Varadhan formula, cf.~Remark~\ref{rem:pinsker},~\ref{rem:DV}, by approximating bounded $\sigma(\omega^{U_\infty})$-measurable functions by $\mathcal F_k$-measurable ones).
\end{proof}

We will apply the above to the loop soup \cite{lawler_werner_2004}, see also the monographs \cite{le_jan_markov_2011,sznitman_topics_2012}, which we now briefly discuss. Let $\rL$ denote the set of (unrooted) loops on $\tilde\cG$, and for $\kappa\geq0$ let $\mu_\kappa$ be the loop measure associated to the diffusion on $\tilde\cG_\kappa$, see \cite[Section~2]{drewitz_cluster_2022}, regarded as a measure on $\rL$ by erasing the excursions of loops into the killing cables, i.e.\ by taking their trace on $\tilde\cG$, cf.\ \cite[(2.31) and below]{drewitz_cluster_2022}. {By \cite[(2.4) and (4.5)]{drewitz_cluster_2022},} the image of $\mu_\kappa$, restricted to loops visiting $G$, under the map associating to a loop its trace on $G$ (which is the continuous-time loop on $G$ obtained by time-changing the loop by the sum of its local times at the vertices of $G$) is the loop measure of the walk on $G$ killed at rate $\kappa$ (recall that $\kappa_x=\kappa\lambda_x$). Parametrising at the speed of $X$, this means that for all measurable $F\geq0$ depending on loops only through their trace on $G$,
\begin{equation}\label{eq:trace}
\int F\,d\mu_\kappa=\sum_{z\in G}\int_0^\infty E_z\big[F\big((X_s)_{0\leq s\leq t}\big);\,X_t=z\big]\,e^{-\kappa t}\,\frac{dt}{t},
\end{equation}
where on the right, $F$ is evaluated at the loop of duration $t$ obtained by identifying the endpoints of $(X_s)_{0\leq s\leq t}$. For $\gamma$ visiting $G$ and $x\in G$, let $\ell_x(\gamma)$ be $\lambda_x^{-1}$ times the time spent at $x$ by the trace of $\gamma$, and for $\kappa:G\to[0,\infty)$ let $\langle\kappa,\ell\rangle=\langle\kappa,\ell(\gamma)\rangle\coloneqq\sum_x\kappa_x\ell_x(\gamma)$. In particular, $\langle\lambda,\ell(\gamma)\rangle$ is the duration of the trace of $\gamma$, and $\langle\kappa,\ell\rangle=\kappa\langle\lambda,\ell\rangle$ when $\kappa_x=\kappa\lambda_x$.

Let $\mathbb Q_\kappa$ denote the law of the loop soup at intensity $\frac12$, i.e.\ of the Poisson point process on $\rL$ with intensity $\frac12\mu_\kappa$, and $\mathbb Q\coloneqq\mathbb Q_0$. Two loops belong to a \emph{cluster} if they can be joined by a finite chain of loops with intersecting ranges. For $x\in G$, let $\mathcal C_x$ denote the union of the ranges of the loops in the cluster containing $x$. Then, by \cite[Proposition~2.1]{lupu_loop_2016}, for all $\kappa\geq0$,
\begin{equation}\label{eq:loop_gff}
\Law_{\mathbb Q_\kappa}(\mathcal C_x)
=
\Law_{\P_\kappa}
\big(
\text{the connected component of }x
\text{ in }\{|\tilde\varphi|>0\}\cap\tilde\cG
\big).
\end{equation}

The following lemma describes how functionals of the duration $\langle\lambda,\ell\rangle$ of a loop going through a \emph{single} point $y \in G$ behave under $\mu_0$. For $y\in G$, let $\rL_y$ be the set of loops in $\rL$ whose range contains $y$. It will be convenient to decompose
\begin{equation}\label{eq:mu0_decomp}
\mu_0=\mu_0^{\nt}+\mu_0^{\triv},
\end{equation}
where $\mu_0^{\nt}$ is the restriction of $\mu_0$ to the loops visiting at least two vertices of $G$.

\begin{lemma}\label{lem:one_point}
For any $y\in G$ and any measurable function $f:[0,\infty)\to[0,\infty)$,
\begin{align}
  \label{eq:general_trivial}
 & \int_{\rL_y} f(\langle\lambda,\ell\rangle) \, d\mu_0^{\triv} = \int_0^\infty \frac{e^{-t}}{t} f(t)\,dt,\quad\text{and}\\
 \label{eq:general_fundamental}
 &\int_{\rL_y} f(\langle\lambda,\ell\rangle) \, d\mu_0^{\nt} \leq \int_0^\infty P_y(X_t=y)\, f(t)\,dt .
\end{align}
\end{lemma}
\begin{proof}
By \eqref{eq:trace}, writing $N_t$ for the number of jumps of $X$ in $[0,t]$ and $(Y_n)_{n\geq0}$ for its discrete skeleton, which is independent of the unit-rate Poisson process $(N_t)_{t\geq0}$, the left-hand side of \eqref{eq:general_trivial} equals
\[
\int_0^{\infty} P_y(N_t=0)f(t)\frac{dt}{t} =  \int_0^\infty \frac{e^{-t}}{t} f(t)\,dt.
\]
For $\mu_0^{\nt}$, we bound $1\{\gamma\in\rL_y\}$ by the number of visits $N_y(\gamma)$ of the trace of $\gamma$ to $y$, which equals $\#\{0\leq i<n:Y_i=y\}$ on $\{N_t=n,\,X_t=z\}$ in \eqref{eq:trace}. By the Markov property, $\sum_zP_z(Y_i=y,\,Y_n=z)=\sum_zP_z(Y_i=y)P_y(Y_{n-i}=z)=P_y(Y_n=y)$ for $0\leq i<n$. Thus, the left-hand side of \eqref{eq:general_fundamental} is bounded by
\begin{align*}
&\sum_{n\geq1}nP_y(Y_n=y)\int_0^\infty P(N_t=n)f(t)\,\frac{dt}{t}=\sum_{n\geq1}P_y(Y_n=y)\int_0^\infty P(N_t=n-1)f(t)\,dt\\
&=\int_0^\infty P_y(Y_{N_t+1}=y)\,f(t)\,dt=\int_0^\infty\Big(\sum_{w\in G}P_y(X_t=w)\,p_{w,y}\Big)\,f(t)\,dt.
\end{align*}
Let $P_tf(x)\coloneqq E_x[f(X_t)]$ and $Pf(x)\coloneqq E_x[f(Y_1)]$ denote the semigroups of $X$ and its discrete skeleton $Y$. Both are self-adjoint on $L^2(\lambda)$ by reversibility, they commute, and $\|P\|_{L^2(\lambda)\to L^2(\lambda)}\leq1$ by Jensen's inequality. Hence, with $g=\lambda_y^{-1}1_{\{y\}}$, the term in parenthesis in the previous display can be recast as $\lambda_y\langle PP_{t/2}g,P_{t/2}g\rangle_{\lambda}\leq\lambda_y\|P_{t/2}g\|^2_{\lambda}=P_y(X_t=y)$. 
\end{proof}

\begin{lemma}\label{lem:loop_density}
For $\kappa\geq0$, regarding $\mu_\kappa$ as a measure on loops of $\tilde\cG$ by erasing their excursions into the killing cables (which does not affect ranges within $\tilde\cG$; see Remark~\ref{rem:sec3_end},\ref{rem:killed_cable_loop} as to why), 
\begin{equation}\label{eq:loop_density}
\frac{d\mu_\kappa}{d\mu_0}=e^{-\langle\kappa,\ell\rangle}.
\end{equation}
\end{lemma}
\begin{proof}
By \eqref{eq:trace}, on loops visiting $G$, the image of $\mu_\kappa$ under the trace map on $G$ has density $e^{-\kappa t}$ with respect to that of $\mu_0$, where $t=\langle\lambda,\ell\rangle$ is the duration of the trace, and $\kappa\langle\lambda,\ell\rangle=\langle\kappa,\ell\rangle$. Moreover, by \cite[Section~2]{drewitz_cluster_2022}, conditionally on its trace on $G$, the loop is obtained by adding Brownian excursions in the cables of $\tilde\cG$ according to its local times at the vertices, with the same law under $\mu_\kappa$ and $\mu_0$. Hence \eqref{eq:loop_density} holds on loops visiting $G$. Loops not visiting $G$ are contained in a single cable of $\tilde\cG$, have the same law under $\mu_\kappa$ and $\mu_0$, and satisfy $\ell=0$.
\end{proof}

\begin{proof}[Proof of Proposition \ref{lem:entropy}]
We now construct an exploration of $\cC_x$ of the form of~\eqref{eq:explore}. 
For any counting measure $\omega$ on $\rL$, we write $\Ran(\omega)=\cup_{\gamma\in \text{supp}(\omega)}\ran(\gamma) (\subseteq\tilde \cG)$ and let $\Ran_x(\omega)$ be the connected component of $x$ in $\Ran(\omega)$.
Now let 
  \begin{align}\label{eq:def_exploration}
  \begin{split}
    &S_0=\big\{\gamma\in \rL:\ran(\gamma)\cap \{x\}\neq \emptyset \big\},\quad \text{and for $n\geq0,$}\\
    &S_{n+1}=\big\{\gamma\in \rL\!\setminus\!  U_{n}: \ran(\gamma)\cap \Ran_x(\omega^{U_{n}}) \neq \emptyset \big\}.
  \end{split}
  \end{align}
In other words, at each step $n\geq 1$, we reveal all previously unexplored loops that intersect the discovered component of $x$. In particular, every loop of the loop cluster of $x$ belongs to $U_\infty$, by induction along a chain of intersecting loops connecting it to a loop through $x$. Note that by construction $\sigma(\mathcal{C}_x)\subseteq \sigma(\omega^{U_\infty})$.

We may assume that $\chi_\kappa(x)<\infty$, for otherwise \eqref{eq:entropy} is trivial.
We first check that $\int_{S_n}(1-e^{-\langle \kappa,\ell\rangle })\,d\mu_0 <\infty$ $\mathbb Q_\kappa$-a.s.\ for all $n\geq0$, as required in Lemma~\ref{lem:poisson}.
Let $\rL_G\coloneqq\bigcup_{y\in G}\rL_y$ denote the set of loops visiting $G$. Loops in $\rL\setminus\rL_G$ satisfy $\ell=0$, hence $e^{-\langle\kappa,\ell\rangle}=1$, and do not contribute to the integrals below. 
Let $\overline{\mathcal{C}}_x\defeq \overline{(\mathcal C_x\cap G)}$, where $\overline D=D\cup\partial D$. By construction, every loop $\gamma\in U_\infty$ (whether or not it belongs to the soup) intersects $\mathcal C_x$. If moreover $\gamma\in U_{\infty}\cap\rL_G$, then $\gamma$ visits a vertex of $\overline{\mathcal{C}}_x$: indeed, either $\gamma$ meets $\mathcal C_x$ at a vertex, which then lies in $\mathcal C_x\cap G$, or only in the interior of a cable $[y,z]$, in which case $\gamma$ passes through $y$ or $z$, and both lie in $\overline{\mathcal C}_x$ since one of them lies in $\mathcal C_x\cap G$. Hence $U_\infty\cap \rL_G \subseteq \bigcup_{y\in \overline{\mathcal{C}}_x} \rL_y$ $\mathbb{Q}_\kappa$-a.s. Furthermore, by~\eqref{eq:loop_gff} and the symmetry of $\tilde\varphi$, and as \eqref{eq:p_0} bounds the degree by $c_0^{-1}$, we have that $E_{\mathbb Q_\kappa}[|\overline{\mathcal{C}}_x|]\leq 2(1+c_0^{-1})\chi_\kappa(x)<\infty$. Consequently, for every measurable $F:\rL\to[0,\infty]$ vanishing on $\rL\setminus\rL_G$, $\mathbb Q_\kappa$-a.s.,
\begin{equation}\label{eq:loop_union_bound}
 \int_{U_\infty}F\,d\mu_0\leq \sum_{y\in \overline{\mathcal C}_x} \int_{\rL_y}F\,d\mu_0.
\end{equation}
Applying~\eqref{eq:general_trivial} and~\eqref{eq:general_fundamental} with $f(t)=1-e^{-\kappa t}$, so that $f(\langle\lambda,\ell\rangle)=1-e^{-\langle\kappa,\ell\rangle}$, gives
\begin{align*}
 &\int_{\rL_y} (1-e^{-\langle\kappa,\ell\rangle}) \, d\mu_0^{\triv} = \int_0^\infty \frac{e^{-t}}{t} (1-e^{-\kappa t})\,dt
\leq \kappa \int_0^\infty e^{-t}\,dt\leq \kappa,\\
&\int_{\rL_y} (1-e^{-\langle\kappa,\ell\rangle}) \, d\mu_0^{\nt} \leq \int_0^\infty P_y(X_t=y)\,dt
= \lambda_yg(y,y)<\infty,
\end{align*}
which, in view of~\eqref{eq:mu0_decomp} and \eqref{eq:loop_union_bound} applied with $F=1-e^{-\langle\kappa,\ell\rangle}$, as $S_n\subseteq U_\infty$ and $|\overline{\mathcal C}_x|<\infty$ $\mathbb Q_\kappa$-a.s., concludes the verification of $\int_{S_n}(1-e^{-\langle \kappa,\ell\rangle })\,d\mu_0 <\infty$, $\mathbb Q_\kappa$-a.s.

We are now ready to apply~\eqref{eq:explore}.
By~\eqref{eq:explore} and \eqref{eq:loop_union_bound} applied with $F=\ophi(e^{-\langle\kappa,\ell\rangle})$ (recall $\ophi(1)=0$),
\[
 2 H\big(\Law_{\mathbb Q_\kappa}(\omega^{U_\infty})\,\big|\,\Law_{\mathbb Q}(\omega^{U_\infty})\big)= E_{\mathbb Q_\kappa}\Big[\int_{U_\infty}\ophi(e^{-\langle\kappa,\ell\rangle})\,d\mu_0\Big]\leq C\chi_\kappa(x)\sup_{y\in G}\int_{\rL_y}\ophi(e^{-\langle\kappa,\ell\rangle})\,d\mu_0.
\]
Let $\sigma$ be an independent random variable taking values in $\{\pm1\}$ each with probability $1/2$.
By~\eqref{eq:loop_gff} and~\eqref{eq:def_exploration}, for every event
$A\in\sigma(\cC_x)$, we can find $A'\in \sigma(\omega^{U_\infty},\sigma)$ such that $\P(A)=\mathbb{Q}(A')$ and $\P_\kappa(A)=\mathbb{Q}_\kappa(A')$. Since $\sigma$ is independent of the loop soup and has the same law under $\mathbb Q$ and $\mathbb Q_\kappa$, the relative entropy of the law of $(\omega^{U_\infty},\sigma)$ under $\mathbb Q_\kappa$ with respect to that under $\mathbb Q$ equals the left-hand side of the above display (e.g.\ by the chain rule, see \cite[Section~5.3]{gray_1990_entro}). Hence
 in view of~\eqref{eq:gen} and the above display, it suffices to show that
\begin{equation}\label{eq:exp_goal}
 \sup_{y\in G}\int_{\rL_y}\ophi(e^{-\langle\kappa,\ell\rangle})\,d\mu_0\leq c \kappa^2 \tilde{\nabla}_\kappa.
\end{equation}
First note that using $1-e^{-x}\leq x$, we have that for all $x\geq 0$,
\begin{equation}\label{eq:ophi}
 \ophi(e^{-x})=1-e^{-x}-xe^{-x}
\leq (x-xe^{-x})\wedge 1=(x(1-e^{-x}))\wedge 1\leq x^2\wedge 1.
\end{equation}
We now apply~\eqref{eq:general_trivial} and~\eqref{eq:general_fundamental} with $f(t)=\ophi(e^{-\kappa t})$, so that $f(\langle\lambda,\ell\rangle)=\ophi(e^{-\langle\kappa,\ell\rangle})$, and bound $f(t)\leq(\kappa t)^2\wedge1=\kappa^2(t\wedge\kappa^{-1})^2$ by~\eqref{eq:ophi}, to obtain for any $y\in G$
\begin{align}
 \label{eq:final_point_nt}
 &\int_{\rL_y}\ophi(e^{-\langle\kappa,\ell\rangle})\,d\mu_0^{\nt}
\leq \kappa^2 \sup_{z\in G}\int_0^\infty P_z(X_t=z)(t\wedge \kappa^{-1})^2\,dt= \kappa^2\tilde{\nabla}_{\kappa},\quad\text{and}\\
\label{eq:final_point_triv}
&\int_{\rL_y}\ophi(e^{-\langle\kappa,\ell\rangle})\,d\mu_0^{\triv} 
\leq \kappa^2 \int_0^{\kappa^{-1}} e^{- t} t \,dt
+ \int_{\kappa^{-1}}^{\infty} e^{- t} \,\frac{dt}{t}
\leq c\kappa^2,
\end{align}
which, in view of~\eqref{eq:mu0_decomp} and since $\tilde\nabla_\kappa\geq\int_0^1t^2e^{-t}\,dt$, concludes the proof of~\eqref{eq:exp_goal} and~\eqref{eq:entropy}.
\end{proof}

\begin{remark}\label{rem:sec3_end}
\begin{enumerate}[label=\arabic*),itemsep=4pt]
\item\label{rem:killed_cable_loop} (Killing cables). Let $L^{\text{kill}}$ be the set of loops $\gamma$ on $\tilde\cG_0$ such that $\gamma$ travels beyond length $(2\kappa_x)^{-1}$ on the open cable attached to $x$ for some $x\in G$.
If one retains the loop excursions in the killing cables in Lemma~\ref{lem:loop_density}, one would instead get $\frac{d\mu_\kappa}{d\mu_0}(\gamma)=1_{\{\gamma\notin L^{\text{kill}}\}}$. As a consequence,~\eqref{eq:final_point_triv} becomes
\[
\int_{\rL_y}\ophi(1_{\{\gamma\notin L^{\text{kill}}\}})\,d\mu_0^{\triv}(\gamma) = \mu_0^{\triv}(\rL_y\cap L^{\text{kill}})=
\int_0^\infty \frac{e^{-t}}{t}(1-e^{-\kappa t})\,dt
=\log(1+\kappa)\asymp \kappa.
\]
After summing over the expected $\asymp\kappa^{-1}$ vertices of $\cC_x\cap G$ as in \eqref{eq:loop_union_bound}, the relative entropy is thus of order $1$ rather than small, so that \eqref{eq:gen} becomes useless.\black
\item\label{rem:field_vs_loops} (Field vs.~loops). It is noteworthy that \eqref{eq:twopoint} (and hence \eqref{eq:vol_bound_general}) are used for $\kappa>0$ and remain pertinent: they yield the volume bound $\chi_\kappa(x)\leq C\kappa^{-1}$ which makes \eqref{eq:entropy} effective, see the proof of Proposition~\ref{prop:ub}. At the level of the excursion sets $\{\tilde\varphi\geq h\}$, the mass shifts the critical point $\tilde h_*$ to a strictly negative level, so the sign cluster becomes subcritical, see \cite[Theorem~4.1]{prevost_counterexamples_2023}, and \cite[Remark~3.4(2)]{rodriguez_zhang_anomalous_2025} for a case in point. This shows that the loop soup is a natural framework here, because the distance from criticality is measured directly in terms of $\kappa$ (cf.~Lemma~\ref{lem:loop_density}), rather than via the shift $|\tilde h_*(\kappa)|$ of the critical level, which is not explicit.\black
\item\label{rem:ber_compa} (Comparison with Bernoulli percolation). For the level perturbation $p<p_c$ used in \cite{dewan_muirhead_2023,vEGPS}, the entropy bound involves the term $(p_c-p)^2\chi(p)$ with no triangle, but the bound $\chi(p)\leq C(p_c-p)^{-1}$ requires the triangle condition \cite{aizenman_newman_1984}, see also \cite{hutchcroft_2022_on} for a recent short proof. Here the roles are reversed: the volume bound is free, while a triangle diagram enters the entropy bound itself, through the second moment of loop durations, see \eqref{eq:final_point_nt}. It seems plausible that a similar effect could be achieved for Bernoulli percolation by perturbing instead with a ghost field \`a la Aizenman--Barsky \cite{aizenman_barsky_1987}, which is a close analogue to the mass.
\end{enumerate}
\end{remark}

\section{A `Simon--Lieb' inequality}\label{sec:SL}

In this section we prove the following: 
\begin{proposition}[`Simon--Lieb' inequality]\label{lem:SL}
Suppose 
\eqref{eq:p_0} holds and  $( \lambda_x \vee g(x,x))\leq C$ for all $x\in G$.
Then for all $x\in G$, finite $D\ni x$, finite $B\subset D^c$ and all $\kappa\in [0,1]$, 
\begin{equation}\label{eq:SL}
{\P_\kappa\big(x\leftrightarrow B\big)\leq C\,E_x\big[e^{-\kappa T_D}\,\P_\kappa\big(X_{T_D}\leftrightarrow B\big)\big].}
\end{equation}
\end{proposition}

The terminology stems from inequalities with a similar flavour originally proved in the context of the Ising model \cite{simon_correlation_1980,lieb_refinement_1980}, see also \cite{duminil-copin_tassion_2016,duminil-copin_alternative_2024,vEGPS,panis_schapira_2026} for recent related developments in the context of Bernoulli percolation.
In the present context the inequality \eqref{eq:SL} will come as a consequence of the Lupu-Werner formula \cite{LW18} for connections and the Markov property of the field. In fact, inequalities of similar flavour have been already derived in \cite[Proposition 3.1]{dpr26fractal} (with obstacles acting as Dirichlet boundary rather than a uniform killing measure) and in the proof of \cite[Proposition 3.1]{cai_ding_2024}.

\begin{remark}
The right-hand side of \eqref{eq:SL} can be expressed in terms of Green functions, as follows. For $D\subset G$ finite, $x,y\in G$ and $\kappa\geq0$, we denote by $g^D_\kappa(x,y)$ the Green function of $X$ killed at rate $\kappa$ and upon exiting $D$ (so that $g^D_\kappa(x,y)=0$ unless $x,y\in D$). It is defined as \eqref{eq:green}, but with the restriction to times $ t \leq T_D\coloneqq \inf\{t\geq0:X_t\notin D\}$. Since $X$ jumps from $y$ to $z$ at
rate $p_{y,z}=\lambda_{y,z}/\lambda_y$ and every $z\in D^c$ with $p_{y,z}>0$ for some $y\in D$
belongs to $\partial D$, 
\begin{equation}\label{eq:exit_id}
\sum_{y\in D}\lambda_{y,z}\,g_\kappa^D(x,y)=\int_0^\infty e^{-\kappa t}\sum_{y\in D}P_x\big(X_t=y,\,t<T_D\big)\,p_{y,z}\,dt
=E_x\big[e^{-\kappa T_D}1\{X_{T_D}=z\}\big],
\end{equation}
and therefore \eqref{eq:SL} can be equivalently recast as
\begin{equation}\label{eq:SL-green}
\P_\kappa\big(x\leftrightarrow B\big)\leq C\sum_{y\in D}\sum_{z\in\partial
D}\lambda_{y,z}\,g_\kappa^D(x,y)\,\P_\kappa\big(z\leftrightarrow B\big);
\end{equation}
this is the form in which \eqref{eq:SL} will be proved. Summing \eqref{eq:exit_id} over $z$ yields $E_x[e^{-\kappa T_D}]$, cf.~\eqref{eq:ub-max}.
For $\kappa=0$ this equals $P_x(T_D<\infty)=1$, so that the weights in \eqref{eq:SL-green} merely average over $\partial D$ and \eqref{eq:SL} yields no gain:
the gain in \eqref{eq:ub} comes from the killing, and Proposition~\ref{lem:entropy} is precisely what
allows to remove it afterwards.
\end{remark}

\begin{proof}
Throughout the proof, we write $P^\kappa_x$, $x\in\tilde\cG_\kappa$, for the law of the diffusion $X$ on $\tilde\cG_\kappa$ starting at $x$, i.e.\ the Brownian motion on $\tilde\cG$ killed upon reaching the open end of one of the killing cables attached at the vertices of $G$, see Section~\ref{sec:prelim}. For a compact set $K\subset\tilde\cG_\kappa$, we write $\tilde\partial K$ for its topological boundary in $\tilde\cG_\kappa$ (denoted by $\partial K$ in \cite{dpr26fractal}, to which we frequently refer below), not to be confused with the vertex boundary $\partial D$ of $D\subset G$.

We follow the arguments in \cite[Section~3]{dpr26fractal} and use the set-to-set connection formula \cite[(18)]{LW18}. Recall $\mathcal C_{K}^{-}$ in \cite[(3.19)]{dpr26fractal}, the set consisting of $K$ and the clusters of $\{ \varphi \leq 0\}$ intersecting it, and set $\mathcal C^-\coloneqq\{x\}\cup\mathcal C^-_B$.
By assumption on $g$, the graph $\cG_\kappa$ satisfies condition (Cap) of \cite{drewitz_cluster_2022}, which in turn yields via \cite[Thm 1.1]{drewitz_cluster_2022} that $\mathcal C^-$ is a.s.~bounded, as required for the next argument to apply. In particular, $\mathcal C^-$ is compatible in the sense of \cite[above (3.16)]{dpr26fractal}.

Using \cite[(3.26)]{dpr26fractal}, applied with $A=\{x\}$ and $A'=B$, the same argument giving (3.27) in the proof of Lemma~3.6 therein, 
combined with the fact that $\{x\leftrightarrow B\}\subset\{\varphi_x>0\}$ and that $x\notin\mathcal C^-_B$ on $\{\varphi_x>0\}$, yields that $\P_\kappa$-a.s.\ on $\{\varphi_x>0\}$, letting $\mathcal{H}\coloneqq\sum_{w \in B \cap \tilde\partial\mathcal C_B^-} h_{\mathcal C^{-}}(x,w) \varphi_w\ (\geq0)$,
\begin{equation}\label{eq:LW-set-2-set}
\P_\kappa\big(x \leftrightarrow B\,\big|\,\mathcal A^+_{\mathcal C^-}\big)=1-e^{-2\varphi_x\mathcal{H}}\leq(1\wedge2\varphi_x\mathcal{H}) 1\{x\notin\mathcal C^-_B\},
\end{equation}
where $\mathcal A^+_{\mathcal C^-}$ is as in \cite[(3.16)]{dpr26fractal} and
\begin{equation}\label{eq:def-eff-cond}
    h_{\mathcal C^-}(x,w) = \lambda^{\kappa,-}_x P^\kappa_x\big(Z_{\tilde H_{\mathcal C^-}}=w,\, \tilde H_{\mathcal C^-}<\infty\big),\qquad w\in G^-; 
\end{equation}
here, as in the proof of \cite[Lemma~3.6]{dpr26fractal}, $G^-=G\cup\tilde\partial\mathcal C^-$ and, under $P^\kappa_x$, $(Z_k)_{k\geq0}$ denotes the chain of successive visits to $G^-$ by the diffusion $X$ and $\tilde H_{\mathcal C^-}\coloneqq\inf\{k\geq1:Z_k\in\mathcal C^-\}$. Thus $Z$ is a random walk on the graph $\mathcal G^-$ with vertex set $G^-$ and weights $\lambda^{\kappa,-}_{\cdot,\cdot}$ (with $\lambda^{\kappa,-}_x=\sum_v\lambda^{\kappa,-}_{x,v}$).

We now integrate out $\varphi_x$ and condition to this effect on the smaller $\sigma$-algebra
$$\mathcal A\coloneqq\mathcal A^+_{\mathcal C^-_B}\ (\subset\mathcal A^+_{\mathcal C^-}),$$
where the inclusion follows from \cite[(3.16)]{dpr26fractal} since $\mathcal C^-_B\subset\mathcal C^-$ and both sets are compatible (alternatively one could also contemplate conditioning on $\varphi_x$ instead; we comment on this below \eqref{eq:cond-final}). Note that, as $\mathcal C^-=\{x\}\cup\mathcal C^-_B$, $\mathcal{H}$ is a function of the set $\mathcal C^-_B$ and of the values of $\varphi$ on $\tilde\partial\mathcal C^-_B$, and thus $\mathcal{H} 1\{x\notin\mathcal C^-_B\}$ is $\mathcal A$-measurable. 
Taking conditional expectations with respect to $\mathcal A$ in the previous display, using that $\mathcal H\geq0$ and Jensen's inequality for the concave function $t\mapsto1\wedge t$, we get that
\begin{equation}\label{eq:cond-A}
\P_\kappa\big(x\leftrightarrow B\,\big|\,\mathcal A\big)\leq\ 1\{x\notin\mathcal C^-_B\}\Big(1\wedge2\mathcal{H}\,\E_\kappa\big[\varphi_x^+\,\big|\,\mathcal A\big]\Big).
\end{equation}
It thus suffices to bound $\mathcal H$ and $\E_\kappa[\varphi_x^+\,|\,\mathcal A]$ separately on $\{x\notin\mathcal C^-_B\}$; we will bound both in terms of the same $\mathcal A$-measurable quantity $\tilde{\mathcal H}$, see \eqref{eq:Lambda-Sigma} and \eqref{eq:cond-mean} below. 

We start with $\mathcal H$, defined above \eqref{eq:LW-set-2-set}. Roughly speaking, we use the Markov property near the boundary of $D$ to decompose $ h_{\mathcal C^{-}}(x,w)$. Since $w \in B \subset D^c$, on the event in \eqref{eq:def-eff-cond}, the random walk must cross the boundary of $D$ via an original edge in $\mathcal G$, and this will make the sums over $y,z$ in \eqref{eq:SL-green} appear.

The first issue is with the weight $\lambda^{\kappa,-}_x$ appearing in \eqref{eq:def-eff-cond}, which may be arbitrarily large when $\tilde\partial\mathcal C^-_B$ gets close to $x$. To deal with this, we decompose according to the first step of $Z$: by the Markov property, the probability on the right of \eqref{eq:def-eff-cond} equals $\sum_vP^\kappa_x(Z_1=v)\,P^\kappa_v(Z_{H_{\mathcal C^-}}=w,\,H_{\mathcal C^-}<\infty)$, where $H_{\mathcal C^-}\coloneqq\inf\{k\geq0:Z_k\in\mathcal C^-\}$ and $v$ ranges over the neighbours of $x$ in the graph $\mathcal G^-$, i.e.\ over the vertices $v\in G^-=G\cup\tilde\partial\mathcal C^-$ with $\lambda^{\kappa,-}_{x,v}>0$. Moreover, for each such $v$, $P^\kappa_v\big(Z_{H_{\mathcal C^-}}=w,\,H_{\mathcal C^-}<\infty\big)=P^\kappa_v\big(X_{H_{\mathcal C^-}}=w,\,H_{\mathcal C^-}<\infty\big)$,
with $H_{\mathcal C^-}$ on the right the hitting time of $\mathcal C^-$ by $X$. Indeed, if $v\in\mathcal C^-$, both probabilities equal $1\{v=w\}$. If $v\notin\mathcal C^-$, the first point of $\mathcal C^-=\{x\}\cup\mathcal C^-_B$ visited by $X$ is either $x$ or a point of $\tilde\partial\mathcal C^-_B$; in both cases it lies in $G^-$, thus it is also the first point of $\mathcal C^-$ visited by $Z$, hence the equality. 

Now, 
for $v\in\tilde\partial\mathcal C^-\setminus\{w\}$ one has $P^\kappa_v(X_{H_{\mathcal C^-}}=w,\,H_{\mathcal C^-}<\infty)=1\{v=w\}=0$, and any other neighbour $v$ of $x$ in $\mathcal{G}^-$ (including $v=w$ when $w\sim x$) belongs to $G$, satisfies $v\sim x$ in $\cG$, and the cable of the edge $\{x,v\}$ contains no point of $\mathcal C^-$ other than $x$ and possibly $v$; in particular this cable is unchanged in $\mathcal G^-$, so that $\lambda^{\kappa,-}_{x,v}=\lambda_{x,v}$, cf.\ \cite[Lemma~2.1]{drewitz_cluster_2022}. Combining this  with the fact that $\lambda^{\kappa,-}_xP^\kappa_x(Z_1=v)\leq\lambda^{\kappa,-}_{x,v}$, we deduce that
\begin{equation}\label{eq:first-step}
h_{\mathcal C^-}(x,w)\leq\sum_{v\in G:\,v\sim x}\lambda_{x,v}\,P^\kappa_v\big(X_{H_{\mathcal C^-_B}}=w,\,H_{\mathcal C^-_B}<\infty\big),
\end{equation}
where we also used that $\mathcal C^-_B\subset\mathcal C^-$. The weight $\lambda^{\kappa,-}_x$ no longer appears in \eqref{eq:first-step}, as desired.

Next, let $T$ be the first time at which $X$ visits a vertex of $G\setminus D$ (with $T=\infty$ if $X$ is killed before). Let $v\in D\cup\partial D$, $w\in B\cap\tilde\partial\mathcal C^-_B$. Under $P_v^\kappa$, on the event $\{X_{H_{\mathcal C^-_B}}=w\}$, since $w\in G\setminus D$, we have $T\leq H_{\mathcal C^-_B}$. Hence, by the strong Markov property at time $T$,
\begin{equation}\label{eq:exit-dec}
P^\kappa_v\big(X_{H_{\mathcal C^-_B}}=w,\,H_{\mathcal C^-_B}<\infty\big)\leq\sum_{z\in\partial D}P^\kappa_v\big(X_T=z,\,T<\infty\big)\,P^\kappa_z\big(X_{H_{\mathcal C^-_B}}=w,\,H_{\mathcal C^-_B}<\infty\big).
\end{equation}
Feeding this into \eqref{eq:first-step}, we claim that
\begin{equation}\label{eq:eff-cond-decomp}
    h_{\mathcal C^{-}}(x,w) \leq 2\lambda_x\sum_{y \in D, z \in \partial D}  \lambda_{y,z}g_\kappa^{D}(x,y)P^\kappa_z\big(X_{H_{\mathcal C_B^-}}=w,\,H_{\mathcal C_B^-}<\infty\big).
\end{equation}
To prove \eqref{eq:eff-cond-decomp}, note that the Markov property at the first step allows one to recombine, for all $z\in\partial D$, (since $x\in D$, the first step of the skeleton occurs no later than $T$, and at time $T$ if it lands in $\partial D$, in which case $T=0$ under $P^\kappa_v$),
\begin{equation*}
\begin{split}
\sum_{v\in G:\,v\sim x}\lambda_{x,v}\,P^\kappa_v\big(X_T=z,\,T<\infty\big)&=(1+\kappa)\lambda_x\,P^\kappa_x\big(X_T=z,\,T<\infty\big);
\end{split}
\end{equation*}
 and \eqref{eq:eff-cond-decomp} now follows by virtue of \eqref{eq:exit_id} upon identifying $P^\kappa_x(X_T=z,\,T<\infty)$ with the expectation on the right-hand side of \eqref{eq:exit_id}.

In view of \eqref{eq:eff-cond-decomp}, we introduce, for $z\in G$,
\begin{equation}\label{eq:eta-tilde}
\tilde\eta(z)\defeq\sum_{w\in B\cap\tilde\partial\mathcal C^-_B}P^\kappa_z\big(X_{H_{\mathcal C^-_B}}=w,\,H_{\mathcal C^-_B}<\infty\big)\,\varphi_w\ (\geq0).
\end{equation}
Recalling the definition of $\mathcal H$ above \eqref{eq:LW-set-2-set}, multiplying \eqref{eq:eff-cond-decomp} by $\varphi_w\geq0$, summing over $w\in B\cap\tilde\partial\mathcal C^-_B$ and using that $\lambda_x\leq C$, we obtain
\begin{equation}\label{eq:Lambda-Sigma}
\mathcal{H}\leq C\tilde{\mathcal H},\qquad\text{where }\tilde{\mathcal H}\defeq\sum_{y\in D}\sum_{z\in\partial D}\lambda_{y,z}\,g^D_\kappa(x,y)\,\tilde\eta(z)\ (\geq0).
\end{equation}
This is the desired bound on $\mathcal{H}$.
We now turn to $\E_\kappa[\varphi_x^+\,|\,\mathcal A]$. By the Markov property \cite[(3.17)]{dpr26fractal} applied to the compatible set $\mathcal C^-_B$, on $\{x\notin\mathcal C^-_B\}$ the conditional law of $\varphi_x$ given $\mathcal A$ is Gaussian with variance at most $g_\kappa(x,x)\leq g(x,x)\leq C$ and mean $\eta^\kappa_{\mathcal C^-_B}(x)$, where $\eta^\kappa_{\mathcal C^-_B}(z)\coloneqq\sum_{y\in\tilde\partial\mathcal C^-_B}P^\kappa_z(X_{H_{\mathcal C^-_B}}=y,\,H_{\mathcal C^-_B}<\infty)\,\varphi_y$, $z\in G$, is the harmonic extension of $\varphi$ on $\mathcal C^-_B$, see \cite[(3.18)]{dpr26fractal}. Since $\varphi$ vanishes on $\tilde\partial\mathcal C^-_B\setminus B$, this mean equals $\tilde\eta(x)$. Arguing as in the proof of \eqref{eq:eff-cond-decomp}, but starting from \eqref{eq:exit-dec} with $v=x$ instead of \eqref{eq:first-step}, we get $0\leq\tilde\eta(x)\leq\tilde{\mathcal H}$. As $(a+b)^+\leq a+b^+$ for $a\geq0$ and $E[\xi^+]\leq\sqrt{\operatorname{Var}(\xi)}$ for a centred Gaussian variable $\xi$, it follows that
\begin{equation}\label{eq:cond-mean}
\E_\kappa\big[\varphi_x^+\,\big|\,\mathcal A\big]\leq\tilde{\mathcal H}+C\qquad\text{on }\{x\notin\mathcal C^-_B\}.
\end{equation}
Plugging \eqref{eq:Lambda-Sigma} and \eqref{eq:cond-mean} into \eqref{eq:cond-A}, and using that $1\wedge C(t^2+t)\leq C't$ for $t\geq0$, we get
\begin{equation}\label{eq:cond-final}
\P_\kappa\big(x\leftrightarrow B\,\big|\,\mathcal A\big)\leq1\wedge C\tilde{\mathcal H}\big(\tilde{\mathcal H}+1\big)\leq C\,\tilde{\mathcal H}.
\end{equation}
Note that both $\mathcal H$ and the conditional mean $\tilde\eta(x)$ of $\varphi_x$ given $\mathcal A$ are bounded by the same quantity $\tilde{\mathcal H}$, see \eqref{eq:Lambda-Sigma} and \eqref{eq:cond-mean}; this is the reason for conditioning on $\mathcal A$ rather than on $\varphi_x$. By contrast, conditioning on $\varphi_x$ shifts the boundary values on $\tilde\partial\mathcal C^-_B$, and hence $\mathcal H$, by an amount which doesn't seem to  be  controllable by $\tilde{\mathcal H}$, and this leads to an additional (additive) term.

 Taking expectations in \eqref{eq:cond-final} yields the desired inequality \eqref{eq:SL}, in the equivalent form \eqref{eq:SL-green}, but with
$\E_\kappa\big[\tilde\eta(z)\big]$ in place of the connection probability on the right-hand side. Thus the desired result is reduced to showing that for all $z\in\partial D$,
\begin{equation}\label{eq:2nd-reduction}
\E_\kappa\big[\tilde\eta(z)\big]\leq C\,\P_\kappa(z\leftrightarrow B).
\end{equation}
By \eqref{eq:eta-tilde}, $\tilde\eta(z)=\eta^\kappa_{\mathcal C^-_B}(z)$ if $z\notin\mathcal C^-_B$, $\tilde\eta(z)=\varphi_z\mathbbm1\{z\in\tilde\partial\mathcal C^-_B\}=\varphi_z^+$ a.s.\ if $z\in B$ (if $\varphi_z>0$, then $\varphi>0$ in a neighbourhood of $z$ in $\tilde\cG_\kappa$, whence $z\in\tilde\partial\mathcal C^-_B$; if $\varphi_z<0$, then $z$ lies in the interior of $\mathcal C^-_B$, and $\tilde\eta(z)=0$ if $z\in\mathcal C^-_B\setminus B$. If $z\in B$, then $\E_\kappa[\tilde\eta(z)]\leq\sqrt{g(z,z)}\leq C\leq2C\,\P_\kappa(\varphi_z\geq0)\leq2C\,\P_\kappa(z\leftrightarrow B)$. If $z\notin B$, then $\tilde\eta(z)=\eta^\kappa_{\mathcal C^-_B}(z)\mathbbm1\{z\notin\mathcal C^-_B\}$, and \eqref{eq:2nd-reduction} is \cite[(3.23)]{dpr26fractal}, applied on the graph $\cG_\kappa$ with $K'=B$ and the compatible set $\mathcal K=\mathcal C^-_B$. Although \cite{dpr26fractal} works under \eqref{eq:Gnu}, the proof of \cite[Lemma~3.8]{dpr26fractal} only uses the Markov property \cite[(3.17)]{dpr26fractal}, the symmetry of $\varphi$ and \cite[Lemma~3.3]{cai_ding_2024}, whose proof (written on $\mathbb Z^d$) applies to any graph with bounded degree, edge weights bounded from above and below and Green function bounded from above and below on the diagonal, cf.\ the end of the proof of \cite[Lemma~3.8]{dpr26fractal}. All of these hold on $\cG_\kappa$, uniformly in $\kappa\in[0,1]$ (in particular $ (2\lambda_x)^{-1} \leq ((1+\kappa)\lambda_x)^{-1}\leq g_\kappa(x,x)\leq g(x,x)\leq C$). 
\end{proof}

\section{Denouement}\label{sec:ub}

In this section we prove Proposition~\ref{prop:ub} and conclude with the proof of Theorem~\ref{thm:main}. 
We start with the former, which follows readily by combining the entropy bound, Proposition~\ref{lem:entropy}, and the Simon--Lieb inequality, Proposition~\ref{lem:SL}, from the previous two sections.

\begin{proof}[Proof of Proposition~\ref{prop:ub}]
Let $A\coloneqq\{x\leftrightarrow B\}$, which belongs to $\sigma(\cC_x)$ and is increasing. By \eqref{eq:vol_bound_general}, as $\sup_{x,y\in G}\tfrac{\lambda_x}{\lambda_y}\leq C$ (see below \eqref{eq:lb}), we have $\chi_\kappa(x)\leq C\kappa^{-1}$. Hence, by
\eqref{eq:entropy} and \eqref{eq:SL},
\begin{equation*}
\begin{split}
\P\big(x\leftrightarrow B\big)&\leq2\,\P_\kappa\big(x\leftrightarrow B\big)+C\kappa\tilde{\nabla}_{\kappa}
\leq2C\,E_x\big[e^{-\kappa T_D}\,\P_\kappa\big(X_{T_D}\leftrightarrow B\big)\big]+C\kappa\tilde{\nabla}_{\kappa},
\end{split}
\end{equation*}
and \eqref{eq:ub} follows since $\kappa\mapsto\P_\kappa(z\leftrightarrow B)$ is non-increasing for every $z\in G$, see \eqref{eq:mono}.
\end{proof}

\begin{proof}[Proof of Theorem~\ref{thm:main}]

The lower bound is an immediate consequence of Corollary~\ref{cor:lb}:
if $2\alpha\leq3\nu$, then $2\beta\leq\nu$, so the right-hand side of \eqref{eq:cor_lb} is
$cR^{-\beta}$.

We now show the upper bound, which will follow from Proposition~\ref{prop:ub} and an induction argument over scales. Recall from \eqref{eq:def-b} that $b(R)\leq C$ if $\epsilon=3\nu-2\alpha>0$ and $b(R)\leq C\log(1+R)$ if $\epsilon=0$. The desired upper bound thus amounts to showing, 
under \eqref{eq:p_0}, \eqref{eq:Valpha} and \eqref{eq:Gnu}, that 
\begin{equation}\label{eq:cor_ub}
\theta(x,R)\leq CR^{-\beta}b(R), \quad \text{$x\in G$, $R\geq1$}.
\end{equation}
Let $\bar{\theta}(R)\coloneqq\sup_{x\in G}\theta(x,R)$, and $R_0\coloneqq A^{1/\beta}\vee4(\CL{C:dist}\vee1)$, where $\CR{C:dist}$ is such that $d\leq\CR{C:dist}d_G$, see \cite[(2.8)]{drewitz_geometry_2024}. Then $B(y,R/2)\subset B(x,R)$ for all $x\in G$, $R\geq R_0$ and $y\in\partial B(x,R/4)$. 
Now let $x\in G$ and $R\geq R_0$, and apply \eqref{eq:ub-max} with $D\coloneqq B(x,R/4)$, $B\coloneqq \partial B(x,R)$ and $\kappa\coloneqq AR^{-\beta}\in(0,1]$. Then $\P(x\leftrightarrow B)=\theta(x,R)$ and, for all $y\in\partial D$, $\P(y\leftrightarrow B)\leq\theta(y,R/2)\leq \bar{\theta}(R/2)$. 
Since $b$ is non-decreasing and $A\geq1$, we have $\kappa\tilde\nabla_\kappa\leq CAR^{-\beta}b(A^{-1/\beta}R)\leq CAR^{-\beta}b(R)$ by \eqref{eq:kappa-tri}, and \eqref{eq:exit} gives $E_x[e^{-\kappa T_D}]\leq Ce^{-cA^{1/\beta}}$. Feeding these bounds into \eqref{eq:ub-max} and taking a sup over $x$ yields
\[
\bar{\theta}(R)\leq Ce^{-cA^{1/\beta}}\bar{\theta}(R/{2})+\CL{C:rec2}AR^{-\beta}b(R),\quad R\geq R_0.
\]
Choose $A$ such that $Ce^{-cA^{1/\beta}}\leq{2}^{-\beta-1}$ and set $\CL{C:ind}\coloneqq2\CR{C:rec2}A\vee(2R_0)^\beta$. For $R\leq2R_0$ we have $\bar\theta(R)\leq1\leq \CR{C:ind}R^{-\beta}b(R)$ using that $b(R) \geq 1$. Now let $R_n=R_0 2^{n+1}$. If $\bar\theta(R_{n-1})\leq \CR{C:ind}R_{n-1}^{-\beta}b(R_{n-1})$ for some $n \geq 1$, then, as $b(R_{n-1})\leq b(R_n)$,
\[
\bar{\theta}(R_n)\leq2^{-\beta-1}\CR{C:ind}2^\beta R_n^{-\beta}b(R_n)+\CR{C:rec2}AR_n^{-\beta}b(R_n)=\big(\tfrac12\CR{C:ind}+\CR{C:rec2}A\big)R_n^{-\beta}b(R_n)\leq \CR{C:ind}R_n^{-\beta}b(R_n).
\]
By induction, \eqref{eq:cor_ub} thus follows, first for $R=R_n$, $n \geq 0$, and all $R \geq 2R_0$ by monotonicity.
\end{proof}

\subsection*{Acknowledgements}
PFR is supported
by the European Research Council (ERC) under the European Union's Horizon Europe research and
innovation programme, grant agreement No.~101171046. JS is supported by an Overseas Bursary from Trinity College and a stipend from DPMMS. WZ is supported by the CDT in Mathematics of Random
Systems and by the President's PhD Scholarship scheme at Imperial College London.

\subsection*{Declaration}
All mathematical results and ideas presented in this article are due to the authors; no artificial intelligence
tools were used in the creative process leading to them. The forthcoming work
\cite{rodriguez_shi_zhang_prep} has benefited from the use of such tools, which will be discussed
separately there.

{\normalsize
\phantomsection\addcontentsline{toc}{section}{References}
\bibliography{bibliography}
\bibliographystyle{abbrv}
}

\end{document}